\documentclass[11pt,reqno,centertags]{amsart}

\usepackage{amsmath,amssymb}
\usepackage{mathrsfs}
\usepackage{tikz-cd}
\usepackage{microtype}
\usepackage[hidelinks]{hyperref}

\makeatletter
\newcommand{\anchorvalue}[2]{%
  \expandafter\gdef\csname paper@number@#1\endcsname{#2}}
\newcommand{\anchor}[1]{%
  \Hy@raisedlink{\hypertarget{#1}{}}%
  \protected@write\@auxout{}{%
    \string\anchorvalue{#1}{\@currentlabel}}}
\DeclareRobustCommand{\anchorref}[1]{%
  \ifcsname paper@number@#1\endcsname
    \hyperlink{#1}{\csname paper@number@#1\endcsname}%
  \else
    \textbf{??}\G@refundefinedtrue
    \@latex@warning{Reference `#1' undefined}%
  \fi}
\DeclareRobustCommand{\eqanchorref}[1]{%
  \textup{\normalfont(\anchorref{#1})}}
\let\paper@pendingeq\@empty
\newcommand{\eqanchor}[1]{%
  \ifmeasuring@\else\gdef\paper@pendingeq{#1}\fi}
\AtBeginDocument{%
  \pretocmd{\tagform@}{%
    \ifmeasuring@\else
      \ifx\paper@pendingeq\@empty\else
        \expandafter\anchor\expandafter{\paper@pendingeq}%
        \global\let\paper@pendingeq\@empty
      \fi
    \fi
  }{}{\PackageError{paper}{Cannot install equation anchors}{}}}
\makeatother

\newcommand{\compactdisplays}{%
  \setlength{\abovedisplayskip}{4pt plus 1pt minus 1pt}%
  \setlength{\belowdisplayskip}{4pt plus 1pt minus 1pt}%
  \setlength{\abovedisplayshortskip}{0pt plus 1pt}%
  \setlength{\belowdisplayshortskip}{2pt plus 1pt}%
  \setlength{\jot}{2pt}%
}

\newtheorem{theorem}{Theorem}
\newtheorem{lemma}[theorem]{Lemma}
\newtheorem{proposition}[theorem]{Proposition}
\theoremstyle{definition}

\newtheorem*{question}{Question}
\newcommand{\M}{\mathcal{M}}
\newcommand{\weak}{L^{p,\infty}(\M,\tau)}
\DeclareMathOperator{\RePart}{Re}

\title[Isometries of noncommutative weak $L^p$]{Surjective linear isometries of noncommutative weak $L^p$}
\author{Ruidong Wang}
\address{School of Mathematical Sciences, Tianjin University of Technology, Tianjin 300384, China}
\email{wangruidong@tjut.edu.cn}
\author[Jianglong Luo]{Jianglong Luo$^{*}$}
\address{School of Mathematical Sciences, Tianjin University of Technology, Tianjin 300384, China}
\email{2100056943@qq.com}
\thanks{$^{*}$Corresponding author: Jianglong Luo
(\href{mailto:2100056943@qq.com}{\nolinkurl{2100056943@qq.com}}).}
\thanks{This work is supported by the Natural Science Foundation of China
(Grant No.~12271402) and the Natural Science Foundation of Tianjin Municipal
Science and Technology Commission (Grant Nos.~23JCYBJC00560 and 22JCYBJC00420).}
\hypersetup{pdfauthor={Ruidong Wang and Jianglong Luo}}

\begin{document}

\begin{abstract}
Isometry classification asks how much of an underlying algebra can be
recovered from the geometry of its function spaces. For noncommutative
$L^p$ spaces, classical results connect surjective linear isometries
with Jordan $*$-isomorphisms. For the natural weak $L^p$ quasi-norm,
Fang, Gao, Huang, and Sukochev recently obtained a classification for
$1<p<\infty$ over $\sigma$-finite semifinite non-atomic von Neumann
algebras, leaving the range $0<p\le1$ open. We settle this remaining
range for $\sigma$-finite non-atomic von Neumann algebras with a faithful
normal semifinite trace $\tau$ satisfying $\tau(1)=\infty$, without
assuming positivity of the isometry. Every surjective complex-linear
isometry has the unique form $T(x)=au\widetilde J(x)$, where $a>0$,
$u$ is unitary, and $J$ is a normal Jordan $*$-automorphism satisfying
$\tau\circ J=a^{-p}\tau$; here $\widetilde J$ denotes its extension
to measurable operators. For finite traces, we obtain
$T(x)=u\widetilde J(x)$ for every $0<p<\infty$, with $J$ trace-preserving
and without any non-atomicity or positivity assumption.
\end{abstract}

\maketitle
\compactdisplays

\section{Introduction}

The classification of linear isometries relates the geometry of a function
space to the measure or algebra from which it is constructed. Weak $L^p$
spaces are natural objects in this problem because of their role in weak-type
estimates and interpolation; see \cite{BS} for the classical setting and
\cite{PX} for noncommutative $L^p$ theory. Here we study isometries for the
natural weak $L^p$ quasi-norm itself.

Weighted composition is a recurring form for surjective linear
isometries of symmetric function spaces. The classical $L^p$ case
for $0<p<\infty$, $p\ne2$, is covered by Lamperti's theorem
\cite[Theorem~3.1]{Lamperti}. For complex symmetric Banach function
spaces whose norm is not proportional to the $L^2$ norm, Zaidenberg
obtained a weighted composition representation
\cite[Theorem~1(b)]{Zaidenberg}. Carothers, Haydon, and Lin classified linear
isometries of weighted Lorentz function spaces $L^{w,p}$ for
$1\le p<\infty$ under their hypotheses on the weight \cite{CHL}.

In the noncommutative setting, Kadison's theorem \cite{Kadison} recovers
the Jordan structure of a unital $C^*$-algebra from its surjective linear
isometries. Yeadon \cite{Yeadon} described isometries of semifinite
noncommutative $L^p$ spaces for $1\le p<\infty$, $p\ne2$.
Sherman \cite{Sherman05} obtained the surjective classification for arbitrary
von Neumann algebras in this range and developed further structure results
for noncommutative $L^p$ isometries in \cite{Sherman06}.

The construction of noncommutative symmetric spaces from singular values
was developed by Dodds, Dodds, and de Pagter \cite{DDP}.
For isometries on these spaces, de Jager and Conradie \cite{dJC20}
obtained representations under positivity and finiteness-preservation
assumptions, and removed positivity under additional assumptions on the
initial space. They also classified surjective isometries of certain
noncommutative Lorentz spaces $L^{w,1}$ and of $L^1+L^\infty$ and
$L^1\cap L^\infty$ \cite{dJC}. Huang and Sukochev treated symmetric
operator spaces with order-continuous, non-Hilbertian norms over
semifinite non-atomic von Neumann algebras or atomic algebras with equal-trace
atoms \cite[Theorem~1.2]{HS}. Huang, Sukochev, and Zanin \cite{HSZ}
treated order-preserving isometries of symmetrically quasi-normed spaces
under strict log-monotonicity assumptions. These results do not directly
cover the natural weak $L^p$ quasi-norm, which is neither order-continuous
nor strictly log-monotone in the non-atomic setting.

For quasi-Banach Marcinkiewicz spaces, Fang, Gao, Huang, and
Sukochev represented positive surjective isometries as scalar
multiples of Jordan $*$-isomorphisms under hypotheses on the
algebras and quasi-norms \cite[Theorem~1.6]{FGHS}.
For the natural weak $L^p$ quasi-norm $Q$ in \eqanchorref{eq:define-Q},
the same authors classified all surjective linear isometries over $\sigma$-finite
semifinite non-atomic von Neumann algebras for $1<p<\infty$ \cite[Theorem~6.9]{FGHS}.
For the commutative space $L^{p,\infty}(0,\infty)$, their result
gives $Tf=a u(f\circ\phi)$, where $|u|=1$ almost everywhere and
the invertible measurable transformation $\phi$ scales Lebesgue
measure by $m(\phi^{-1}(A))=a^{-p}m(A)$
\cite[Corollary~6.11]{FGHS}.

In \cite[Remark~6.10]{FGHS}, Fang, Gao, Huang, and Sukochev leave open
the classification for $0<p\le1$ over non-atomic von Neumann algebras. The problem
can be stated as follows.

\begin{question}
For $0<p\le1$, characterize the surjective linear isometries, including
the positive ones, of $L^{p,\infty}(\M,\tau)$ when $\M$ is a non-atomic von Neumann algebra.
\end{question}

The following theorem and the finite-trace result in
Section~\anchorref{sec:finite-trace} answer this question under the algebra
and trace assumptions stated below, without assuming that the isometry
is positive. We first treat $\sigma$-finite
non-atomic von Neumann algebras with $\tau(1)=\infty$.

\begin{theorem}\anchor{thm:main}
Let $0<p\le 1$ and let $\M$ be a $\sigma$-finite non-atomic von Neumann algebra with faithful normal semifinite trace $\tau$ satisfying $\tau(1)=\infty$. A complex-linear map $T\colon\weak\to\weak$ is a surjective isometry for the natural quasi-norm if and only if
\[
T(x)=a u\widetilde J(x),\qquad x\in\weak,
\]
for a scalar $a>0$, a unitary $u\in\M$, and a normal Jordan $*$-automorphism $J$ of $\M$ satisfying $\tau\circ J=a^{-p}\tau$. Here $\widetilde J$ is the extension to $\tau$-measurable operators. The triple $(a,u,J)$ is unique.
\end{theorem}

The finite-trace case is proved separately in Section~\anchorref{sec:finite-trace},
for all $0<p<\infty$ and without assuming that the algebra is non-atomic.

For Lebesgue measure on $(0,\infty)$, Theorem~\anchorref{thm:main}
gives the same weighted composition formula as
\cite[Corollary~6.11]{FGHS} for $0<p\le1$.

\section{Preliminaries}

\subsection{Projections and measurable operators}

We use the standard theory of noncommutative integration developed by
Segal \cite{Segal} and Nelson \cite{Nelson}, with the notation for
measurable operators in \cite[Section~1]{FK}. The algebra $\M$ acts on
$\mathcal H$, and $\mathcal P(\M)$ denotes its projections. For projections,
$e\le f$ means $ef=fe=e$, and $e\vee f$ denotes their least upper bound.
The corner $e\M e$ has identity $e$.

The finite-trace projections form a directed set and increase strongly
to $1$. We use
\[
\tau(e\vee f)\le\tau(e)+\tau(f)
\]
and the standard trace-splitting property of a non-atomic von Neumann algebra:
\[
e\in\mathcal P(\M),\quad 0\le\lambda<\tau(e)
\quad\Longrightarrow\quad
\exists f<e\text{ a projection with }\tau(f)=\lambda.
\]

Write $S(\M,\tau)$ for the $*$-algebra of $\tau$-measurable operators
affiliated with $\M$, with sums and products understood as their closures
\cite[Section~1]{FK}. For $x=v|x|\in S(\M,\tau)$, its right and left
supports are
\[
r(x)=v^*v=1_{(0,\infty)}(|x|),\qquad
\ell(x)=vv^*=r(x^*),\qquad \tau(r(x))=\tau(\ell(x)).
\]
Here $v\in\M$ is a partial isometry, and $1_I(h)$ denotes the spectral
projection of a self-adjoint affiliated operator $h$ for a Borel set $I$.
We write $\RePart x=(x+x^*)/2$. Orthogonality means
$x^*y=xy^*=0$, equivalently $r(x)r(y)=\ell(x)\ell(y)=0$.

As in \cite[Section~2.1]{FGHS}, the bounded finite-support operators form
the $*$-subalgebra
\[
\mathcal F=\{x\in\M:\tau(r(x))<\infty\}
=\bigcup_{\substack{e\in\mathcal P(\M)\\\tau(e)<\infty}}e\M e.
\]
For $x\in\mathcal F$, one may take $e=r(x)\vee\ell(x)$, since
$exe=x$ and $\tau(e)\le2\tau(r(x))<\infty$. Any two finite-trace
corners are contained in $(e\vee f)\M(e\vee f)$.

\subsection{Singular values and weak \texorpdfstring{$L^p$}{Lp}}

For $x\in S(\M,\tau)$, define
\begin{equation}\eqanchor{eq:count-values}
d_x(\lambda)=\tau(1_{(\lambda,\infty)}(|x|)),\qquad
\mu_s(x)=\inf\{\lambda\ge0:d_x(\lambda)\le s\}.
\end{equation}
Here $\lambda\ge0$ and $s>0$, and $d_x(0)$ may be infinite.
Both functions are nonincreasing and right-continuous. We use
\cite[Proposition~2.2 and Lemma~2.5]{FK}
\[
d_x(\mu_s(x))\le s,\qquad
d_x(\lambda)>s\quad(0\le\lambda<\mu_s(x)).
\]
For $0<p<\infty$, the natural weak $L^p$ quasi-norm is
\begin{equation}\eqanchor{eq:define-Q}
\weak=\{x\in S(\M,\tau):Q(x)<\infty\},\qquad
Q(x)=\sup_{s>0}s^{1/p}\mu_s(x).
\end{equation}
This is a quasi-normed space \cite[Section~1.2]{FGHS}. Throughout,
an isometry is a complex-linear map preserving $Q$.

We use the projection formula and singular value inequalities
\cite[Definition~2.1, Proposition~2.2 and Lemma~2.5]{FK}:
\begin{equation}\eqanchor{eq:compute-mu}
\mu_s(x)=\inf\{\|xe\|_\infty:e\in\mathcal P(\M),\ \tau(1-e)\le s\},
\end{equation}
where unbounded cuts have norm $+\infty$, and
\begin{align}
\mu_{s+t}(x+y)&\le\mu_s(x)+\mu_t(y),\eqanchor{eq:sum}\\
\mu_s(bxc)&\le\|b\|_\infty\|c\|_\infty\mu_s(x)
\qquad(b,c\in\M),\eqanchor{eq:multiply-bound}\\
\mu_s(x^*)&=\mu_s(|x|)=\mu_s(x).\eqanchor{eq:star-values}
\end{align}
In particular, for unitary $u,v$, scalars $\alpha$, and bounded $b$,
\[
\mu_s(uxv)=\mu_s(x),\qquad
Q(\alpha x)=|\alpha|Q(x),\qquad
\mu_s(x+b)\le\mu_s(x)+\|b\|_\infty.
\]
The distribution formula is \cite[Eq.~(2.4)]{BCLJ}
\begin{equation}\eqanchor{eq:compute-Q}
Q(x)=\sup_{\lambda>0}\lambda\,d_x(\lambda)^{1/p}.
\end{equation}
For $|x|=ce$, with $c>0$ and $e$ a finite-trace projection, this gives
\[
\mu_s(x)=c\mathbf1_{(0,\tau(e))}(s),\qquad Q(x)=c\tau(e)^{1/p}.
\]
We also write $\|x\|_j^j=\tau(|x|^j)$ for $j=2,4$. These moments are
finite on $\mathcal F$, and when $|x|=ce$,
\[
\|x\|_j^j=c^j\tau(e)=c^{j-p}Q(x)^p.
\]

\subsection{Convergence in measure and Jordan maps}

We only need the following criteria for convergence in measure
\cite[Lemma~3.1]{FK}:
\begin{align*}
x_n\to x\text{ in measure}
&\quad\Longleftrightarrow\quad
\mu_s(x_n-x)\to0\text{ for every }s>0\\
&\quad\Longleftrightarrow\quad
d_{x_n-x}(\epsilon)\to0\text{ for every }\epsilon>0.
\end{align*}
Also, $Q$ is lower semicontinuous in measure \cite[Lemma~3.4]{FK}:
\[
x_n\to x\text{ in measure}\quad\Longrightarrow\quad
Q(x)\le\liminf_nQ(x_n).
\]

A Jordan $*$-automorphism is a complex-linear bijection preserving squares
and adjoints. Its Jordan product is $x\circ y=(xy+yx)/2$, and it is normal
if it preserves suprema of bounded increasing nets of positive elements.
It scales the trace by $\kappa>0$ if
$\tau(J(x))=\kappa\tau(x)$ for $x\in\M_+$, allowing infinite values.

Such maps preserve order, projections, and orthogonality of projections.
We use the identity, obtained from
$exe=2e\circ(e\circ x)-e\circ x$:
\begin{equation}\eqanchor{eq:compress-J}
J(exe)=J(e)J(x)J(e).
\end{equation}
If $J$ is normal and scales the trace by $\kappa$, then
\cite[Lemma~2.2]{dJC}, applied with target trace $\kappa^{-1}\tau$,
gives its unique measure-continuous extension $\widetilde J$ to
$S(\M,\tau)$. The homomorphism and antihomomorphism decomposition
\cite{Stormer}, together with $\mu_s(x^*)=\mu_s(x)$, gives
\begin{equation}\eqanchor{eq:scale-Q}
d_{\widetilde J(x)}(\lambda)=\kappa d_x(\lambda)\quad(\lambda>0),
\qquad Q(\widetilde J(x))=\kappa^{1/p}Q(x).
\end{equation}
Consequently, if $\kappa=a^{-p}$ and $u$ is unitary, then
$V(x)=au\widetilde J(x)$ satisfies
\[
Q(V(x))=a\kappa^{1/p}Q(x)=Q(x).
\]
Its inverse is $V^{-1}(y)=\widetilde J^{-1}(a^{-1}u^*y)$, which also
preserves $Q$ by \eqanchorref{eq:scale-Q}. Thus $V$ is a surjective
complex-linear isometry of $\weak$.

\section{Proof of the main theorem}

Throughout this section, assume that $\tau(1)=\infty$ and $0<p\le1$.
For $x\in\weak$ and $t>0$, write
$\mu_{t-}(x)=\lim_{s\uparrow t}\mu_s(x)$. Define
\[
D_0(x)=\limsup_{s\downarrow0}s^{1/p}\mu_s(x),\qquad
D_\infty(x)=\limsup_{s\to\infty}s^{1/p}\mu_s(x).
\]
\begin{lemma}\anchor{lem:attain-Q}
For every $x\in\weak$,
if $\max\{D_0(x),D_\infty(x)\}<Q(x)$, then there is some
$t\in(0,\infty)$ with $Q(x)=t^{1/p}\mu_{t-}(x)$.
\end{lemma}

\begin{proof}
Assume $\max\{D_0(x),D_\infty(x)\}<Q(x)$. By the definition of $Q$,
choose $s_n>0$ with $s_n^{1/p}\mu_{s_n}(x)\to Q(x)$.
If a subsequence $(s_{n_k})$ converges to $0$, then
\[
Q(x)=\lim_{k\to\infty}s_{n_k}^{1/p}\mu_{s_{n_k}}(x)\le D_0(x).
\]
If a subsequence $(s_{n_k})$ tends to $\infty$, then
\[
Q(x)=\lim_{k\to\infty}s_{n_k}^{1/p}\mu_{s_{n_k}}(x)\le D_\infty(x).
\]
Both conclusions contradict the strict inequality. We assume $s_n\to t\in(0,\infty)$.

Fix $0<s<t$. For all sufficiently large $n$, we have $s_n>s$.
By monotonicity, $\mu_{s_n}(x)\le\mu_s(x)$. Hence
\(
s_n^{1/p}\mu_{s_n}(x)\le s_n^{1/p}\mu_s(x).
\)
Letting $n\to\infty$ gives $Q(x)\le t^{1/p}\mu_s(x)$.
Then letting $s\uparrow t$ gives
\[
Q(x)\le t^{1/p}\mu_{t-}(x).
\]
Conversely, by the definition of $Q$, we have $s^{1/p}\mu_s(x)\le Q(x)$
for every $0<s<t$. Taking the left limit at $t$ gives
\[
t^{1/p}\mu_{t-}(x)\le Q(x).
\]
Combining the two inequalities gives $t^{1/p}\mu_{t-}(x)=Q(x)$.
\end{proof}

\begin{lemma}\anchor{lem:unchanged-ends}
For every $x\in\weak$ and $b\in\mathcal F$,
\begin{equation}\eqanchor{eq:unchanged-ends}
D_i(x+b)=D_i(x),\qquad i\in\{0,\infty\}.
\end{equation}
\end{lemma}

\begin{proof}
To prove the identity for $D_0$, fix $s>0$. For every projection
$e$ with $\tau(1-e)\le s$ for which $xe$ is bounded, we have
\[
\|(x+b)e\|_\infty
\le \|xe\|_\infty+\|be\|_\infty
\le \|xe\|_\infty+\|b\|_\infty.
\]
Taking the infimum over $e$ in \eqanchorref{eq:compute-mu} and applying
the same estimate to $x=(x+b)+(-b)$ gives
\begin{align*}
\mu_s(x+b)&\le\mu_s(x)+\|b\|_\infty,\\
\mu_s(x)&\le\mu_s(x+b)+\|b\|_\infty.
\end{align*}
Thus
\(
\bigl|s^{1/p}\mu_s(x+b)-s^{1/p}\mu_s(x)\bigr|
\le s^{1/p}\|b\|_\infty\rightarrow0, 
\) $s\rightarrow0$.
Therefore $D_0(x+b)=D_0(x)$.

To prove the identity for $D_\infty$, recall that $\tau(r(b))<\infty$.
If $\tau(r(b))=0$, faithfulness gives $r(b)=0$, hence $b=0$.
Suppose $\tau(r(b))>0$. Since $\tau(1-(1-r(b)))=\tau(r(b))$,
from \eqanchorref{eq:compute-mu} we have
\(
0\le\mu_{\tau(r(b))}(b)\le\|b(1-r(b))\|_\infty=0.
\)
Hence, for $s>\tau(r(b))$, \eqanchorref{eq:sum} gives
\[
\mu_s(x+b)\le\mu_{s-\tau(r(b))}(x)+\mu_{\tau(r(b))}(b)=\mu_{s-\tau(r(b))}(x).
\]
Multiplying by $s^{1/p}$ we hold
\begin{align*}
s^{1/p}\mu_s(x+b)
 &\le s^{1/p}\mu_{s-\tau(r(b))}(x)\\
 &=\left(\frac{s}{s-\tau(r(b))}\right)^{1/p}
 (s-\tau(r(b)))^{1/p}\mu_{s-\tau(r(b))}(x)\\
 &\leq\left(\frac{s}{s-\tau(r(b))}\right)^{1/p}Q(x), 
\end{align*}
let $s\rightarrow\infty$, we get $D_{\infty}(x+b)\leq Q(x)$.
Applying the same estimate to $x=(x+b)+(-b)$ gives
\[
D_\infty(x)\le D_\infty(x+b).
\]
Thus $D_\infty(x+b)=D_\infty(x)$.
\end{proof}

\begin{lemma}\anchor{lem:distance}
For $x\in \weak$, let $x_n=x1_{[1/n,n]}(|x|)$. Then
\[
\inf_{b\in\mathcal F}Q(x-b)=\lim_{n\to\infty}Q(x-x_n)
   =\max\{D_0(x),D_\infty(x)\}.
\]
\end{lemma}

\begin{proof}
For every $b\in\mathcal F$, Lemma~\anchorref{lem:unchanged-ends} gives
\[
Q(x-b)\ge\max\{D_0(x-b),D_\infty(x-b)\}
        =\max\{D_0(x),D_\infty(x)\}.
\]
The $x_n$ satisfies $\|x_n\|_\infty\le n$ and
\[
\tau(r(x_n))\le d_x(1/(2n))\le(2nQ(x))^p<\infty,
\]
so $x_n\in\mathcal F$. Moreover,
\begin{align*}
\mu_s(x-x_n)
 &=\mu_s\bigl(x(1-1_{[1/n,n]}(|x|))\bigr)\\
 &\le\|1-1_{[1/n,n]}(|x|)\|_\infty\,\mu_s(x)
 \le\mu_s(x).
\end{align*}

To prove the upper bound, suppose that
\[
\max\{D_0(x),D_\infty(x)\}<c<\limsup_nQ(x-x_n).
\]
By definition of $Q(x-x_n)$, we can choose
$s_n>0$ such that $c<s_n^{1/p}\mu_{s_n}(x-x_n)\leq s_n^{1/p}\mu_{s_n}(x)$. Then we can take $\lambda_n$ such that $c<s_n^{1/p}\lambda_n\leq s_n^{1/p}\mu_{s_n}(x-x_n)\leq s_n^{1/p}\mu_{s_n}(x)$.
Passing to a further subsequence, one of the following two cases holds.
If $\lambda_n\ge1/n$, we have
\begin{align*}
0<s_n<d_{x-x_n}(\lambda_n)
 &=d_x(\max\{n,\lambda_n\})\le d_x(n)\le(Q(x)/n)^p\longrightarrow0.
\end{align*}
Hence $c\le\limsup_n s_n^{1/p}\mu_{s_n}(x)
\le D_0(x)$, a contradiction. If $\lambda_n<1/n$, from $c<s_n^{1/p}\lambda_n$ we know 
\[
s_n^{1/p}>\frac{c}{\lambda_n}>cn\longrightarrow\infty,
\]
so $c\le\limsup_n s_n^{1/p}\mu_{s_n}(x)
\le D_\infty(x)$, again a contradiction. Together with the lower bound, this yields
\begin{align*}
\max\{D_0(x),D_\infty(x)\}
 &\le \inf_{b\in\mathcal F}Q(x-b)\\
 &\le \liminf_{n\to\infty}Q(x-x_n)
 \le \limsup_{n\to\infty}Q(x-x_n)\\
 &\le \max\{D_0(x),D_\infty(x)\}.
\end{align*}
\end{proof}

For $r>0$ we define $P_r(x)=\sup_{Q(z)\le r}Q(x+z)$, where
$z\in\weak$. If $T$ is a surjective linear isometry, then
$T$ maps $\{z\in\weak:Q(z)\le r\}$ onto itself. Changing variables
by $T$ and using linearity and preservation of $Q$ gives
\(
P_r(Tx)
=\sup_{Q(z)\le r}Q(Tx+Tz)
=\sup_{Q(z)\le r}Q(T(x+z))
=\sup_{Q(z)\le r}Q(x+z)
=P_r(x).
\)

\begin{lemma}\anchor{lem:largest-sum}
For $x\in \weak$ and $r>0$, we have
\begin{equation}\eqanchor{eq:largest-sum}
P_r(x)=\sup_{t>0}\inf_{0<s<t}t^{1/p}
       \left(\mu_s(x)+\frac{r}{(t-s)^{1/p}}\right).
\end{equation}
\end{lemma}

\begin{proof}
\emph{Upper bound.}
Fix $z$ with $Q(z)\le r$. For every $t>0$ and $0<s<t$,
\eqanchorref{eq:sum} and \eqanchorref{eq:define-Q} give
\begin{align*}
t^{1/p}\mu_t(x+z)
 &\le t^{1/p}\bigl(\mu_s(x)+\mu_{t-s}(z)\bigr)\\
 &\le t^{1/p}\left(\mu_s(x)+\frac{r}{(t-s)^{1/p}}\right).
\end{align*}
The left-hand side is independent of $s$. Hence
\[
Q(x+z)
\le \sup_{t>0}\inf_{0<s<t}
 t^{1/p}\left(\mu_s(x)+\frac{r}{(t-s)^{1/p}}\right).
\]
Taking the supremum over $z$ with $Q(z)\le r$ gives the upper bound
for $P_r(x)$.

\emph{Lower bound.}
Fix $t>0$ and choose
\[
0<d<\inf_{0<s<t}\left(\mu_s(x)+\frac{r}{(t-s)^{1/p}}\right).
\]
We construct $z$ with $Q(z)\le r$ and $Q(x+z)\ge dt^{1/p}$.
The definition of $\mu_t(x)$ gives
\[
\tau(1_{(\mu_t(x),\infty)}(|x|))\le t
\le\tau(1_{[\mu_t(x),\infty)}(|x|)).
\]
Since $\M$ is a non-atomic von Neumann algebra, choose a projection $e$ commuting with $|x|$ such that
\[
1_{(\mu_t(x),\infty)}(|x|)\le e
\le1_{[\mu_t(x),\infty)}(|x|),\qquad\tau(e)=t.
\]
Write $x=v|x|$. If $\mu_t(x)>0$, then $e\le v^*v$.
If $\mu_t(x)=0$, then $\tau(r(x))\le t<\infty$; extend $v$ to a
unitary by \cite[Lemma~2.5]{dJC}. Thus in both cases $v$ acts
isometrically on $r(x)\vee e$. Set
\[
z=v(d-|x|)_+e,\qquad |z|=(d-|x|)_+e.
\]
To estimate $Q(z)$, fix $0<u<t$ and
$\ell>(d-\mu_{t-u}(x))_+$. If $\ell<d$, then
$d_x(d-\ell)>t-u$, and the choice of $e$ gives
\begin{align*}
d_z(\ell)
 &=\tau\bigl(e1_{[0,d-\ell)}(|x|)\bigr)\\
 &\le t-\tau\bigl(e1_{(d-\ell,\infty)}(|x|)\bigr)\\
 &=t-\min\{t,d_x(d-\ell)\}<u.
\end{align*}
If $\ell\ge d$, then $d_z(\ell)=0$. Hence the definition of
$\mu_u$ and the choice of $d$ imply
\[
\mu_u(z)\le(d-\mu_{t-u}(x))_+\le r u^{-1/p}
\qquad(0<u<t).
\]
Since $r(z)\le e$, we also have $\mu_u(z)=0$ for $u\ge t$.
Therefore $Q(z)\le r$.

On the other hand,
\[
|x+z|=|x|(1-e)+\max\{|x|,d1\}e\ge de.
\]
Since $\tau(e)=t$, \eqanchorref{eq:compute-Q} gives
$P_r(x)\ge Q(x+z)\ge\sup_{0<s<t}s^{1/p}\mu_{s}(x+z)\ge\sup_{0<s<t}s^{1/p}d=dt^{1/p}$. Letting $d$ increase to the
chosen infimum yields
\[
P_r(x)\ge t^{1/p}\inf_{0<s<t}
\left(\mu_s(x)+\frac{r}{(t-s)^{1/p}}\right).
\]
Taking the supremum over $t>0$ proves the lower bound.
\end{proof}

\begin{proposition}\anchor{prop:compare-sums}
Let $x,y\in\weak$. If $\mu_s(x)\le\mu_s(y)$ for every $s>0$, then
$P_r(x)\le P_r(y)$ for every $r>0$.
\end{proposition}

\begin{proof}
Fix $r>0$ and $t>0$. For every $0<s<t$, the hypothesis gives
\[
t^{1/p}\left(\mu_s(x)+\frac{r}{(t-s)^{1/p}}\right)
\le
t^{1/p}\left(\mu_s(y)+\frac{r}{(t-s)^{1/p}}\right).
\]
Taking the infimum over $0<s<t$ on both sides gives
\[
\inf_{0<s<t}t^{1/p}\left(\mu_s(x)+\frac{r}{(t-s)^{1/p}}\right)
\le
\inf_{0<s<t}t^{1/p}\left(\mu_s(y)+\frac{r}{(t-s)^{1/p}}\right).
\]
This holds for every $t>0$. Taking the supremum over $t>0$ yields $P_r(x)\le P_r(y)$.
\end{proof}

Next, we define
\(
\mathscr V=\bigl\{x\in\mathcal F: Q(x)=1,\ 
  \exists\,c>0,\ e\in\mathcal P(\M)\ \text{with } |x|=ce\bigr\}.
\)
For $x\in\mathscr V$ with $|x|=ce$, we have $e\ne0$ since $Q(x)=1$.
By spectral calculus, $r(x)=e$ and $\|x\|_\infty=c$, it follows that
\[
|x|=\|x\|_\infty r(x)\qquad(x\in\mathscr V).
\]

\begin{proposition}\anchor{prop:equal-sums}
For every $y\in\mathscr V$ and $r>0$,
\begin{equation}\eqanchor{eq:minimum}
P_r(y)=\max_{t>1}\min\left\{t^{1/p}+r,\ r\left(\frac{t}{t-1}\right)^{1/p}\right\},
\end{equation}
this means any $y\in\mathscr V$, and fixed $r>0$, then $P_r(y)$ is a constant.
\end{proposition}

\begin{proof}
For $y\in\mathscr V$ with $|y|=ce$, the equality
$1=Q(y)=c\tau(e)^{1/p}$ gives
\[
\mu_s(y)=\tau(e)^{-1/p}\mathbf1_{(0,\tau(e))}(s).
\]
Rescaling $s$ and $t$ by $\tau(e)$ in \eqanchorref{eq:largest-sum} yields
\[
P_r(y)=\sup_{t>0}\inf_{0<s<t}
\left(t^{1/p}\mathbf1_{(0,1)}(s)
      +r\left(\frac{t}{t-s}\right)^{1/p}\right),
\]
which is independent of $y$.

For fixed $t$, $t/(t-s)$ increases with $s$.
If $0<t\le1$, the inner infimum is $t^{1/p}+r$, obtained as
$s\downarrow0$. If $t>1$, splitting at $s=1$ gives
\begin{align*}
\inf_{0<s<1}\left(t^{1/p}+r\left(\frac{t}{t-s}\right)^{1/p}\right)
 &=t^{1/p}+r,\\
\inf_{1\le s<t}r\left(\frac{t}{t-s}\right)^{1/p}
 &=r\left(\frac{t}{t-1}\right)^{1/p}.
\end{align*}
Thus the inner infimum is the minimum of these two values.
On $t>1$, the first increases continuously from $1+r$ to $\infty$,
while the second decreases continuously from $\infty$ to $r$.
Their minimum therefore attains its maximum at their unique intersection
$t_*>1$, where
\[
t_*^{1/p}+r=r\left(\frac{t_*}{t_*-1}\right)^{1/p}>1+r.
\]
This exceeds the supremum $1+r$ over $0<t\le1$, so only $t>1$
contributes to the maximum. This proves \eqanchorref{eq:minimum}.
\end{proof}

\begin{lemma}\anchor{lem:recognize-form}
Let $Q(x)=1$ and $y\in\mathscr V$.\\(1) $P_r(x)\ge P_r(y)$ for every
$r>0$. \\(2) If $Q(x)=b^{1/p}\mu_{b-}(x)$ for some $b\in(0,\infty)$
and $P_r(x)=P_r(y)$ for some $r>0$, then $x\in\mathscr V$.
\end{lemma}

\begin{proof}
(1)
Fix $r>0$ and $b>0$ with $\mu_{b-}(x)>0$. Since $\M$ is a non-atomic von Neumann algebra, using
$\tau(1)=\infty$, choose a projection $e$ with $\tau(e)=b$.
The operator $\mu_{b-}(x)e$ is positive and vanishes on $(1-e)\mathcal H$.
Thus
\[
1_{(\lambda,\infty)}\bigl(\mu_{b-}(x)e\bigr)
=
\begin{cases}
e, & 0\le\lambda<\mu_{b-}(x),\\
0, & \lambda\ge\mu_{b-}(x).
\end{cases}
\]
Taking the trace gives $d_{\mu_{b-}(x)e}(\lambda)=b$ in the
first case and $0$ in the second. Therefore \eqanchorref{eq:count-values}
gives
\begin{align*}
\mu_s\bigl(\mu_{b-}(x)e\bigr)
 &=\inf\{\lambda\ge0:d_{\mu_{b-}(x)e}(\lambda)\le s\}\\
 &=
\begin{cases}
\inf[\mu_{b-}(x),\infty)=\mu_{b-}(x), & 0<s<b,\\
\inf[0,\infty)=0, & s\ge b.
\end{cases}
\end{align*}
Hence by the monotonicity of $\mu_s$ we have
\[
\mu_s(\mu_{b-}(x)e)
=\mu_{b-}(x)\mathbf{1}_{(0,b)}(s)\le\mu_s(x)
\qquad(s>0).
\]
 By Proposition~\anchorref{prop:compare-sums},
$P_r(x)\ge P_r(\mu_{b-}(x)e)$. Also,
\(b^{1/p}\mu_{b-}(x)=\lim_{s\uparrow b}s^{1/p}\mu_s(x)\le1.
\)
Thus from $Q(z)\le r$ we have $Q(b^{1/p}\mu_{b-}(x)z)=b^{1/p}\mu_{b-}(x)Q(z)\le r$.
Consequently,
\begin{align*}
P_r(x)
 &\ge P_r(\mu_{b-}(x)e)\\
 &\ge \sup_{Q(z)\le r} Q\bigl(\mu_{b-}(x)e+b^{1/p}\mu_{b-}(x)z\bigr)\\
 &=b^{1/p}\mu_{b-}(x)P_r(b^{-1/p}e)\\
 &=b^{1/p}\mu_{b-}(x)P_r(y).
\end{align*}
The last equality follows from Proposition~\anchorref{prop:equal-sums}
because $b^{-1/p}e\in\mathscr V$. Since $Q(x)=1$, choose $b_n>0$
with $b_n^{1/p}\mu_{b_n}(x)\to1$. Then we have $b_n^{1/p}\mu_{b_n}(x)\le b_n^{1/p}\mu_{b_n-}(x)\le1$, so $b_n^{1/p}\mu_{b_n-}(x)\to1$ as $n\to\infty$.

Taking $b=b_n$ in the lower bound and letting $n\to\infty$
gives $P_r(x)\ge P_r(y)$.

(2)
Suppose $b^{1/p}\mu_{b-}(x)=1$ and $P_r(x)=P_r(y)$ for the
chosen $b,r>0$. By monotonicity, for $0<s<1$, we have $b^{1/p}\mu_{bs}(x)\ge1$. 
Rescaling both variables in \eqanchorref{eq:largest-sum} by $b$ gives
\[
P_r(x)=\sup_{t>0}\inf_{0<s<t}\left((bt)^{1/p}\mu_{bs}(x)+r\left(\frac{t}{t-s}\right)^{1/p}\right).
\]
By Proposition~\anchorref{prop:equal-sums}, there is a unique $t_*>1$
such that
\(
t_*^{1/p}+r=r\left(\frac{t_*}{t_*-1}\right)^{1/p}=P_r(y).
\)
Suppose that $b^{1/p}\mu_{b\delta}(x)>1$ for some $0<\delta<1$, and we 
set $\eta=b^{1/p}\mu_{b\delta}(x)-1>0$. Since $s\mapsto\mu_s(x)$ is
nonincreasing, for $0<s<\delta$, we have $b^{1/p}\mu_{bs}(x)\ge b^{1/p}\mu_{b\delta}(x)=1+\eta$.
For $t>1$ this gives
\begin{align*}
&(bt)^{1/p}\mu_{bs}(x)+r\left(\frac{t}{t-s}\right)^{1/p}
\ge
\begin{cases}
t^{1/p}(1+\eta)+r, & 0<s<\delta,\\
t^{1/p}+r\left(\dfrac{t}{t-\delta}\right)^{1/p}, & \delta\le s<1.
\end{cases}
\end{align*}
At $t=t_*$, both bounds strictly exceed $t_*^{1/p}+r=P_r(y)$. So there is a $t\in(1,t_*)$ for which
\(
\inf_{0<s<1}\left((bt)^{1/p}\mu_{bs}(x)+r\left(\frac{t}{t-s}\right)^{1/p}\right)
>P_r(y).
\)
For the same $t$, since $t/(t-1)$ decreases with $t$, we have
\[
\inf_{1\le s<t}\left((bt)^{1/p}\mu_{bs}(x)+r\left(\frac{t}{t-s}\right)^{1/p}\right)
\ge r\left(\frac{t}{t-1}\right)^{1/p}>P_r(y).
\]
Taking the minimum of these infima gives a value strictly larger
than $P_r(y)$, contradicting $P_r(x)=P_r(y)$.
Together with $b^{1/p}\mu_{bs}(x)\ge1$ for $0<s<1$,
we have $b^{1/p}\mu_{bs}(x)=1$ for all $0<s<1$.

Next suppose that $\mu_{bs}(x)>0$ for some $s>1$. By monotonicity,
there exist $\eta>0$ and $0<\delta<t_*-1$ such that
$b^{1/p}\mu_{bs}(x)\ge\eta$ for $1\le s\le1+\delta$.
For $t>1+\delta$ this gives
\begin{align*}
&(bt)^{1/p}\mu_{bs}(x)+r\left(\frac{t}{t-s}\right)^{1/p}
\ge
\begin{cases}
t^{1/p}\eta+r\left(\dfrac{t}{t-1}\right)^{1/p}, & 1\le s\le1+\delta,\\[3pt]
r\left(\dfrac{t}{t-1-\delta}\right)^{1/p}, & 1+\delta<s<t.
\end{cases}
\end{align*}
At $t=t_*$, both bounds exceed $r(t_* /(t_*-1))^{1/p}=P_r(y)$.
Thus there is a $t>t_*$ for which
\[
\inf_{1\le s<t}\left((bt)^{1/p}\mu_{bs}(x)+r\left(\frac{t}{t-s}\right)^{1/p}\right)
>P_r(y).
\]
For the same $t$, and $b^{1/p}\mu_{bs}(x)\ge1$ for $0<s<1$
gives
\[
\inf_{0<s<1}\left((bt)^{1/p}\mu_{bs}(x)+r\left(\frac{t}{t-s}\right)^{1/p}\right)
\ge t^{1/p}+r>P_r(y).
\]
Again the minimum of the two infima is strictly larger than
$P_r(y)$, contradicting $P_r(x)=P_r(y)$. Therefore $\mu_{bs}(x)=0$ for
every $s>1$. By right continuity,
\[
\mu_b(x)=\lim_{s\downarrow1}\mu_{bs}(x)=0.
\]

Combining these facts with $\mu_{bs}(x)=b^{-1/p}$ for $0<s<1$ gives
\[
\mu_s(x)=
\begin{cases}
b^{-1/p}, & 0<s<b,\\
0, & s\ge b.
\end{cases}
\]
By \eqanchorref{eq:count-values} and right continuity of $d_x$,
$\mu_s(x)>\lambda$ if and only if $s<d_x(\lambda)$. Hence
\[
d_x(\lambda)=
\begin{cases}
b, & 0<\lambda<b^{-1/p},\\
0, & \lambda\ge b^{-1/p}.
\end{cases}
\]
In particular, $\tau(r(x))=\lim_{\lambda\downarrow0}d_x(\lambda)=b$.
For every $0<\lambda<b^{-1/p}$,
\[
\tau\bigl(r(x)-1_{(\lambda,\infty)}(|x|)\bigr)=b-b=0.
\]
By faithfulness, $1_{(\lambda,\infty)}(|x|)=r(x)$ for all such
$\lambda$, and $1_{(b^{-1/p},\infty)}(|x|)=0$.
Letting $\lambda\uparrow b^{-1/p}$ therefore gives
$1_{\{b^{-1/p}\}}(|x|)=r(x)$. Thus $|x|=b^{-1/p}r(x)$,
so $x\in\mathcal F$ and, since $Q(x)=1$, we get $x\in\mathscr V$.
\end{proof}

\begin{lemma}\anchor{lem:split}
Suppose $Q(x)=1$ and $\max\{D_0(x),D_\infty(x)\}=1$. There are
$x_1,x_2\in\weak$ with $x=x_1+x_2$, $Q(x_1)=Q(x_2)=1$, and
\[
\mu(x_1+e^{i\theta}x_2)=\mu(x)\quad(\theta\in\mathbb R),\qquad
\mu(x_j)\le\mu(x)\quad(j=1,2).
\]
\end{lemma}

\begin{proof}
If $D_0(x)=1$, choose $s_n\downarrow0$ with $s_n^{1/p}\mu_{s_n}(x)\to1$.
Otherwise $D_\infty(x)=1$, so choose $s_n\to\infty$ with the same
limit $s_n^{1/p}\mu_{s_n}(x)\to1$. For all sufficiently large $n$,
we have $\mu_{s_n}(x)>0$. For such $n$, set
$\lambda_n=(1-1/n)\mu_{s_n}(x)$. Since $0<\lambda_n<\mu_{s_n}(x)$,
\eqanchorref{eq:count-values} gives $d_x(\lambda_n)>s_n$. Hence, by
\eqanchorref{eq:compute-Q},
\[
(1-1/n)^p\bigl(s_n^{1/p}\mu_{s_n}(x)\bigr)^p
=\lambda_n^p s_n
<\lambda_n^p d_x(\lambda_n)\le Q(x)^p=1.
\]
Therefore $\lambda_n^p d_x(\lambda_n)\to1$.

If $s_n\to\infty$, then
\(
0<\lambda_n\le s_n^{-1/p}\rightarrow0\), and 
\(d_x(\lambda_n)>s_n\rightarrow\infty.
\)
Thus we can pass to a subsequence, with
$\lambda_n<\lambda_{n-1}$ and $d_x(\lambda_n)>n d_x(\lambda_{n-1})$
for $n\ge2$. For these $n$, set
$e_n=1_{(\lambda_n,\lambda_{n-1}]}(|x|)$. Then
\begin{align*}
\tau(e_n)&=d_x(\lambda_n)-d_x(\lambda_{n-1}),\\
\lambda_n^p\tau(e_n)
 &=\lambda_n^pd_x(\lambda_n)
 \left(1-\frac{d_x(\lambda_{n-1})}{d_x(\lambda_n)}\right).
\end{align*}
Since $0\le d_x(\lambda_{n-1})/d_x(\lambda_n)<1/n$, it follows that
$\lambda_n^p\tau(e_n)\to1$.

If $s_n\downarrow0$, then $s_n^{1/p}\lambda_n\to1$, so
$\lambda_n\to\infty$. Moreover,
\[
0<d_x(\lambda_n)\le\lambda_n^{-p}\longrightarrow0.
\]
Thus we can pass to a subsequence, relabeled by $n$, with
$\lambda_{n+1}>\lambda_n$ and $d_x(\lambda_{n+1})<d_x(\lambda_n)/n$
for $n\ge2$. For these $n$, set
$e_n=1_{(\lambda_n,\lambda_{n+1}]}(|x|)$. Then
\begin{align*}
\tau(e_n)&=d_x(\lambda_n)-d_x(\lambda_{n+1}),\\
\lambda_n^p\tau(e_n)
 &=\lambda_n^pd_x(\lambda_n)
 \left(1-\frac{d_x(\lambda_{n+1})}{d_x(\lambda_n)}\right).
\end{align*}
Since $0\le d_x(\lambda_{n+1})/d_x(\lambda_n)<1/n$, again
$\lambda_n^p\tau(e_n)\to1$. In either case the intervals defining
the $e_n$ are pairwise disjoint, so $e_ne_m=0$ for $n\ne m$.

Let $q_1=\sum_{n\ {\rm odd}}e_n$, where the sum is the strong sum
of orthogonal projections, and set $q_2=1-q_1$, $x_j=xq_j$.
Since $q_1+q_2=1$, we have $x_1+x_2=x$. Both $q_j$ commute with
$|x|$, so $|x_j|=|x|q_j$. If $e_n\le q_j$, then
\[
1_{(\lambda_n,\infty)}(|x_j|)
=1_{(\lambda_n,\infty)}(|x|)q_j\ge e_n.
\]
Taking traces and applying \eqanchorref{eq:compute-Q} gives
\[
Q(x_j)^p\ge\lambda_n^p d_{x_j}(\lambda_n)
\ge\lambda_n^p\tau(e_n).
\]
By the definition of $q_1$ and the orthogonality of the $e_n$, we have $e_{2n+1}\le q_1$ and $e_{2n}\le q_2$ for $n\ge1$.
Hence
\begin{align*}
Q(x_1)^p&\ge\lim_{n\to\infty}\lambda_{2n+1}^p\tau(e_{2n+1})=1,\\
Q(x_2)^p&\ge\lim_{n\to\infty}\lambda_{2n}^p\tau(e_{2n})=1.
\end{align*}
Moreover,
\eqanchorref{eq:multiply-bound} gives $\mu_s(x_j)\le\mu_s(x)$ for every $s>0$,
so $Q(x_j)\le Q(x)=1$. Therefore $Q(x_1)=Q(x_2)=1$. 
Since
$q_j^*=q_j=q_j^2$ and $q_1q_2=q_2q_1=0$, we have \((q_1+e^{i\theta}q_2)^*(q_1+e^{i\theta}q_2)=(q_1+e^{i\theta}q_2)(q_1+e^{i\theta}q_2)^*=1\)
Thus $q_1+e^{i\theta}q_2$ is unitary. Since right multiplication by a
unitary preserves singular values, we obtain
\(
\mu_s(x_1+e^{i\theta}x_2)
=\mu_s\bigl(x(q_1+e^{i\theta}q_2)\bigr)
=\mu_s(x)
\) for $s>0$, $\theta\in\mathbb R$.
\end{proof}

\begin{lemma}\anchor{lem:constant-norm}
Let $A,B\in\mathcal F$. Suppose that for every $\theta\in\mathbb R$,
\(
0\ne|A+e^{i\theta}B|
=\|A+e^{i\theta}B\|_\infty r(A+e^{i\theta}B).
\)
If $Q(A+e^{i\theta}B)=Q(A+B)$ for every $\theta\in\mathbb R$, then \(
\|A+e^{i\theta}B\|_\infty=\|A+B\|_\infty
\) for $\theta\in\mathbb R
$
\end{lemma}

\begin{proof}
By hypothesis, $|A+e^{i\theta}B|=\|A+e^{i\theta}B\|_\infty
r(A+e^{i\theta}B)$. Since $r(A+e^{i\theta}B)$ is a projection,
for $j=2,4$ we have
\begin{align*}
\|A+e^{i\theta}B\|_j^j
 &=\tau\bigl(|A+e^{i\theta}B|^j\bigr)
 =\|A+e^{i\theta}B\|_\infty^j\,\tau(r(A+e^{i\theta}B)).
\end{align*}
Since $Q(x)=c\tau(e)^{1/p}$ for $|x|=ce$ and the hypothesis on $Q$,
we also have
\[
Q(A+B)^p=Q(A+e^{i\theta}B)^p
=\|A+e^{i\theta}B\|_\infty^p\,\tau(r(A+e^{i\theta}B)).
\]
Eliminating the trace gives
\[
\|A+e^{i\theta}B\|_j^j
=Q(A+B)^p\|A+e^{i\theta}B\|_\infty^{j-p}
\qquad(j=2,4).
\]
In particular,
\begin{equation}\eqanchor{eq:norm-from-square}
\|A+e^{i\theta}B\|_2^2
=Q(A+B)^p\|A+e^{i\theta}B\|_\infty^{2-p}.
\end{equation}
Since $A+B\ne0$, we have $Q(A+B)>0$. Set $\gamma=(4-p)/(2-p)>2$, we have
\begin{equation}\eqanchor{eq:fourth-and-second}
\begin{aligned}
\|A+e^{i\theta}B\|_4^4
 &=Q(A+B)^p\bigl(\|A+e^{i\theta}B\|_\infty^{2-p}\bigr)^\gamma\\
 &=Q(A+B)^p
   \left(\frac{\|A+e^{i\theta}B\|_2^2}{Q(A+B)^p}\right)^\gamma\\
 &= Q(A+B)^{p(1-\gamma)}
                    \bigl(\|A+e^{i\theta}B\|_2^2\bigr)^\gamma.
\end{aligned}
\end{equation}
and
\begin{equation}\eqanchor{eq:expand-square}
\begin{aligned}
\|A+e^{i\theta}B\|_2^2
 &=\tau((A+e^{i\theta}B)^*(A+e^{i\theta}B)\\
 &=\tau(A^*A+B^*B+e^{i\theta}A^*B+e^{-i\theta}B^*A)\\
 &=\tau(A^*A+B^*B)+2\RePart\bigl(e^{i\theta}\tau(A^*B)\bigr).
\end{aligned}
\end{equation}
Squaring the product gives
\begin{align*}
&\bigl((A+e^{i\theta}B)^*(A+e^{i\theta}B)\bigr)^2\\
&\quad=(A^*A+B^*B)^2+A^*BB^*A+B^*AA^*B\\
&\qquad\quad+e^{i\theta}
 \bigl((A^*A+B^*B)A^*B+A^*B(A^*A+B^*B)\bigr)\\
&\qquad\quad+e^{-i\theta}
 \bigl((A^*A+B^*B)B^*A+B^*A(A^*A+B^*B)\bigr)\\
&\qquad\quad+e^{2i\theta}(A^*B)^2+e^{-2i\theta}(B^*A)^2.
\end{align*}
Taking the trace gives
\begin{align*}
\|A+e^{i\theta}B\|_4^4
 =\tau\!\left(\bigl((A+e^{i\theta}B)^*(A+e^{i\theta}B)\bigr)^2\right)=\sum_{k=-2}^{2} a_k e^{ik\theta},
\end{align*}
where the coefficients $a_k$ are independent of $\theta$.
Since $\|A+e^{i\theta}B\|_4^4$ is real, $a_0\in\mathbb R$ and
$a_{-k}=\overline{a_k}$ for $k=1,2$, we hold
\begin{align*}
\|A+e^{i\theta}B\|_4^4
 &=a_0+2\operatorname{Re}(a_1)\cos\theta
       -2\operatorname{Im}(a_1)\sin\theta\\
 &\quad+2\operatorname{Re}(a_2)\cos2\theta
       -2\operatorname{Im}(a_2)\sin2\theta.
\end{align*}

Suppose $\|A+e^{i\theta}B\|_2^2$ is not constant. If $\tau(A^*B)=0$,
then \eqanchorref{eq:expand-square} gives
\[
\|A+e^{i\theta}B\|_2^2=\tau(A^*A+B^*B)
\qquad(\theta\in\mathbb R),
\]
which is independent of $\theta$. Hence $\tau(A^*B)\ne0$. 
As $\theta$ ranges over $\mathbb R$, we may assume $\tau(A^*B)>0$, and \eqanchorref{eq:expand-square} gives
\[
\|A+e^{i\theta}B\|_2^2
=\tau(A^*A+B^*B)+2\tau(A^*B)\cos\theta.
\]

Since $\cos(-\theta)=\cos\theta$, \eqanchorref{eq:fourth-and-second} gives
$\|A+e^{i\theta}B\|_4^4=\|A+e^{-i\theta}B\|_4^4$.
Hence the sine coefficients vanish,
so there are real constants $\alpha_0,\alpha_1,
\alpha_2$ such that
\begin{align*}
\|A+e^{i\theta}B\|_4^4&=
Q(A+B)^{p(1-\gamma)}
\bigl(\tau(A^*A+B^*B)+2\tau(A^*B)\cos\theta\bigr)^\gamma\\
&=\alpha_0+\alpha_1\cos\theta+\alpha_2\cos(2\theta)
\end{align*}
Putting $t=\cos\theta$ we have
\[
\bigl(\tau(A^*A+B^*B)+2\tau(A^*B)t\bigr)^\gamma
=\frac{\alpha_0+\alpha_1t+\alpha_2(2t^2-1)}{Q(A+B)^{p(1-\gamma)}}
\qquad(-1\le t\le1).
\]
For $-1<t<1$, since $\tau(A^*B)>0$, we have
\begin{align*}
\tau(A^*A+B^*B)+2\tau(A^*B)t
 &>\tau(A^*A+B^*B)-2\tau(A^*B)\\
 &=\|A-B\|_2^2>0.
\end{align*}
Applying the chain rule to the left-hand side of the equality gives
\begin{align*}
&\frac{\mathrm d}{\mathrm dt}
 \bigl(\tau(A^*A+B^*B)+2\tau(A^*B)t\bigr)^\gamma\\
&\qquad=\gamma\bigl(2\tau(A^*B)\bigr)
 \bigl(\tau(A^*A+B^*B)+2\tau(A^*B)t\bigr)^{\gamma-1},\\
&\frac{\mathrm d^2}{\mathrm dt^2}
 \bigl(\tau(A^*A+B^*B)+2\tau(A^*B)t\bigr)^\gamma\\
&\qquad=\gamma(\gamma-1)\bigl(2\tau(A^*B)\bigr)^2
 \bigl(\tau(A^*A+B^*B)+2\tau(A^*B)t\bigr)^{\gamma-2},\\
&\frac{\mathrm d^3}{\mathrm dt^3}
 \bigl(\tau(A^*A+B^*B)+2\tau(A^*B)t\bigr)^\gamma\\
&\qquad=\gamma(\gamma-1)(\gamma-2)\bigl(2\tau(A^*B)\bigr)^3
 \bigl(\tau(A^*A+B^*B)+2\tau(A^*B)t\bigr)^{\gamma-3}>0,
\end{align*}
where the strict inequality uses $\gamma>2$ and the positivity
just proved. Since the right-hand side is a polynomial of degree
at most two in $t$, its third derivative is zero. 
This contradiction proves that $\|A+e^{i\theta}B\|_2^2$ is independent
of $\theta$, so it equals its value at $\theta=0$. Applying
\eqanchorref{eq:norm-from-square} at $\theta=0$ gives
\begin{align*}
Q(A+B)^p\|A+e^{i\theta}B\|_\infty^{2-p}
 &=\|A+e^{i\theta}B\|_2^2\\
 &=\|A+B\|_2^2\\
 &=Q(A+B)^p\|A+B\|_\infty^{2-p}.
\end{align*}
This is $\|A+e^{i\theta}B\|_\infty=\|A+B\|_\infty$. We complete the proof.
\end{proof}

\begin{lemma}\anchor{lem:keep-form}
Every surjective complex-linear isometry of $\weak$ maps $\mathscr V$ onto itself.
\end{lemma}

\begin{proof}
This proof is a bit long, we divide it into several steps.
\textit{Step 1. Construct $A,B\in\mathscr V$.}
Let $e\in\mathscr V$. Since $P_r(T(e))=P_r(e)$, 
by Lemma~\anchorref{lem:recognize-form}, it suffices to prove
\(
\exists\,t\in(0,\infty)\ \text{with }t^{1/p}\mu_{t-}(T(e))=1.
\)
Suppose $t^{1/p}\mu_{t-}(T(e))<1$ for every $t>0$. Then
Lemma~\anchorref{lem:attain-Q} gives
\[
\max\{D_0(T(e)),D_\infty(T(e))\}=1.
\]
By Lemma~\anchorref{lem:split}, we have $T(e)=x_1+x_2$ with $Q(x_j)=1$ and $\mu_s(x_j)\le\mu_s(T(e))$ for $j=1,2$. Furthermore, Lemma~\anchorref{lem:recognize-form} and Proposition~\anchorref{prop:compare-sums}
give 
\(
P_r(e)\le P_r(x_j)\le P_r(T(e))=P_r(e)
\) for $j=1,2,\ r>0.$
Hence $P_r(x_j)=P_r(e)$.
We set $A=T^{-1}(x_1)$ and $B=T^{-1}(x_2)$. Then $A+B=e$,
$Q(A)=Q(B)=1$, and $P_r(A)=P_r(B)=P_r(e)$.
Taking $\theta=\pi$ in Lemma~\anchorref{lem:split} gives
\(
\mu_s(T(e))=\mu_s(x_1-x_2)
\) for $s>0$.
Hence, by \eqanchorref{eq:define-Q} and the isometry of $T$,
\[
Q(A-B)=Q(x_1-x_2)
=\sup_{s>0}s^{1/p}\mu_s(T(e))=Q(T(e))=1.
\]
Since $e\in\mathcal F$, from Lemma~\anchorref{lem:unchanged-ends} and homogeneity of $D_i$ give
\begin{align*}
2D_i(A)&=D_i(2A-e)=D_i(A-B)\le1,\\
2D_i(B)&=D_i(e-2B)=D_i(A-B)\le1
\end{align*}
for $i\in\{0,\infty\}$. For each $z\in\{A,B\}$, we have
$\max\{D_0(z),D_\infty(z)\}\le\tfrac12<Q(z)=1$.
Lemma~\anchorref{lem:attain-Q} therefore gives
\[
\exists\,t\in(0,\infty)\ \text{with }t^{1/p}\mu_{t-}(z)=1.
\]
Lemma~\anchorref{lem:recognize-form} now gives $A,B\in\mathscr V\subset\mathcal F$.

\smallskip\noindent\textit{Step 2. Prove $C_\theta\in\mathscr V$ for every $\theta$.}
Thus for every $\theta$, $C_\theta=A+e^{i\theta}B$ belongs to $\mathcal F$.
By complex linearity,
\[
T(C_\theta)=T(A)+e^{i\theta}T(B)=x_1+e^{i\theta}x_2.
\]
Hence Lemma~\anchorref{lem:split} gives
\[
\mu_s(T(C_\theta))=\mu_s(x_1+e^{i\theta}x_2)=\mu_s(T(e))
\qquad(s>0).
\]
By \eqanchorref{eq:define-Q}, and isometry of $T$, we have
\begin{align*}
Q(C_\theta)&=Q(T(C_\theta))=Q(T(e))=1,\\
P_r(C_\theta)&=P_r(T(C_\theta))=P_r(T(e))=P_r(e)
\qquad(r>0).
\end{align*}
Using the projection $1$ in \eqanchorref{eq:compute-mu} gives
$\mu_s(C_\theta)\le\|C_\theta\|_\infty$ for every $s>0$. Thus
\[
0\le D_0(C_\theta)
\le\lim_{s\downarrow0}s^{1/p}\|C_\theta\|_\infty=0.
\]
For $s\ge\tau(r(C_\theta))$, using $1-r(C_\theta)$ in
\eqanchorref{eq:compute-mu} gives
\[
0\le\mu_s(C_\theta)\le
\bigl\|C_\theta(1-r(C_\theta))\bigr\|_\infty=0.
\]
Since $\tau(r(C_\theta))<\infty$, we obtain
\[
D_\infty(C_\theta)=\limsup_{s\to\infty}s^{1/p}\mu_s(C_\theta)=0.
\]
Therefore $\max\{D_0(C_\theta),D_\infty(C_\theta)\}=0<Q(C_\theta)=1$.
By Lemma~\anchorref{lem:attain-Q},
\[
\exists\,t\in(0,\infty)\ \text{with }t^{1/p}\mu_{t-}(C_\theta)=1.
\]
Lemma~\anchorref{lem:recognize-form} yields $C_\theta\in\mathscr V$.

\smallskip\noindent\textit{Step 3. Compare the norm identities and lead a contradiction.}
Since $C_\theta\in\mathscr V$ and $Q(C_\theta)=1$ for every $\theta$,
Lemma~\anchorref{lem:constant-norm} gives $\|C_\theta\|_\infty=\|A+B\|_\infty$.
By \eqanchorref{eq:expand-square},
\[
\|C_\theta\|_2^2=\|A\|_2^2+\|B\|_2^2
 +2\RePart(e^{i\theta}\tau(A^*B))
\qquad(\theta\in\mathbb R).
\]
Integrating the identity we have
\begin{align*}
\frac{1}{2\pi}\int_0^{2\pi}\|C_\theta\|_2^2\,d\theta
&=\frac{1}{2\pi}\int_0^{2\pi}(\|A\|_2^2+\|B\|_2^2)\,d\theta+\frac{1}{\pi}\int_0^{2\pi}
       \RePart\bigl(e^{i\theta}\tau(A^*B)\bigr)\,d\theta\\
&=\|A\|_2^2+\|B\|_2^2.
\end{align*}
For $A,B,C_\theta\in\mathscr V$, we have
\begin{align*}
\|A\|_2^2&=Q(A)^p\|A\|_\infty^{2-p}=\|A\|_\infty^{2-p},\\
\|B\|_2^2&=Q(B)^p\|B\|_\infty^{2-p}=\|B\|_\infty^{2-p},\\
\|C_\theta\|_2^2&=Q(C_\theta)^p\|A+B\|_\infty^{2-p}=\|A+B\|_\infty^{2-p}.
\end{align*}
Substituting these into the integral identity yields
\begin{equation}\eqanchor{eq:strict-bound}
\|A+B\|_\infty^{2-p}=\|A\|_\infty^{2-p}+\|B\|_\infty^{2-p}.
\end{equation}

Since $|C_\theta|=\|A+B\|_\infty r(C_\theta)$, we have
\begin{equation}\eqanchor{eq:expand-cube}
C_\theta C_\theta^*C_\theta
=\|A+B\|_\infty^2 C_\theta
=\|A+B\|_\infty^2(A+e^{i\theta}B).
\end{equation}
On the other hand, 
\begin{align*}
		C_\theta C_\theta^*C_\theta
		&=e^{-i\theta}AB^*A
		+(AA^*A+AB^*B+BB^*A)\\
		&\quad+e^{i\theta}(AA^*B+BA^*A+BB^*B)
		+e^{2i\theta}BA^*B.
\end{align*}
Multiplying \eqanchorref{eq:expand-cube} by $e^{-2i\theta}$ and
expanding $C_\theta C_\theta^*C_\theta$ gives
\begin{equation}\eqanchor{eq:weight-cube}
\begin{aligned}
e^{-2i\theta}C_\theta C_\theta^*C_\theta
 &=e^{-3i\theta}AB^*A
       +e^{-2i\theta}(AA^*A+BB^*A+AB^*B)\\
 &\quad+e^{-i\theta}(AA^*B+BA^*A+BB^*B)+BA^*B\\
 &=\|A+B\|_\infty^2(e^{-2i\theta}A+e^{-i\theta}B).
\end{aligned}
\end{equation}
Integrating \eqanchorref{eq:weight-cube} over $[0,2\pi]$ and
dividing by $2\pi$ gives
\begin{align*}
BA^*B
 &=\frac{\|A+B\|_\infty^2}{2\pi}
   \int_0^{2\pi}(e^{-2i\theta}A+e^{-i\theta}B)\,d\theta=0.
\end{align*}
Integrating \eqanchorref{eq:expand-cube} over $[0,2\pi]$ and dividing
by $2\pi$ gives
\begin{equation}\eqanchor{eq:constant-terms}
\begin{aligned}
AA^*A+BB^*A+AB^*B
 =\frac{\|A+B\|_\infty^2}{2\pi}
   \int_0^{2\pi}(A+e^{i\theta}B)\,d\theta
 =\|A+B\|_\infty^2 A.
\end{aligned}
\end{equation}
Since $A\in\mathscr V$, we have $|A|=\|A\|_\infty r(A)$.
Both sides of \eqanchorref{eq:constant-terms}
gives
\begin{equation}\eqanchor{eq:mixed-terms}
\begin{aligned}
BB^*A+AB^*B
 &=\|A+B\|_\infty^2 A-AA^*A\\
 &=\|A+B\|_\infty^2 A-\|A\|_\infty^2 A\\
 &=(\|A+B\|_\infty^2-\|A\|_\infty^2)A.
\end{aligned}
\end{equation}

We now use \eqanchorref{eq:mixed-terms} to prove
$\|A+B\|_\infty^2=\|A\|_\infty^2+\|B\|_\infty^2$.
Set
\[
\delta=\|A+B\|_\infty^2-\|A\|_\infty^2>0,
\]
where positivity follows from \eqanchorref{eq:strict-bound}. Then
\eqanchorref{eq:mixed-terms} becomes
\[
\delta A=BB^*A+AB^*B.
\]
Since $B\in\mathscr V$, its squared moduli satisfy
\[
(BB^*)^2=\|B\|_\infty^2BB^*,\qquad
(B^*B)^2=\|B\|_\infty^2B^*B.
\]
Moreover, taking adjoints in $BA^*B=0$ gives $B^*AB^*=0$, so
\[
BB^*AB^*B=B(B^*AB^*)B=0.
\]
Multiplying the identity for $\delta A$ on the left by $BB^*$ and
on the right by $B^*B$, and adding, therefore gives
\begin{align*}
\delta^2A
 &=BB^*(\delta A)+(\delta A)B^*B\\
 &=(BB^*)^2A+2BB^*AB^*B+A(B^*B)^2\\
 &=\|B\|_\infty^2(BB^*A+AB^*B)
 =\|B\|_\infty^2\delta A.
\end{align*}
As $\delta>0$ and $A\ne0$, this implies $\delta=\|B\|_\infty^2$, hence
\[
\|A+B\|_\infty^2=\|A\|_\infty^2+\|B\|_\infty^2.
\]
Finally, $0<1-p/2<1$ and $A,B\ne0$ give
\[
\|A+B\|_\infty^{2-p}
=\bigl(\|A\|_\infty^2+\|B\|_\infty^2\bigr)^{1-p/2}
<\|A\|_\infty^{2-p}+\|B\|_\infty^{2-p},
\]
contradicting \eqanchorref{eq:strict-bound}. Thus $T(e)\in\mathscr V$.
Applying the same argument to $T^{-1}$ shows that
$T(\mathscr V)=\mathscr V$.
\end{proof}

Let $T$ be a surjective complex-linear isometry of
$\weak$. Our next goal is to find $a>0$, a unitary $u\in\M$, and
a normal Jordan $*$-automorphism $J$ of $\M$ such that
\[
T(x)=auJ(x)\qquad(x\in\mathcal F),\qquad
\tau\circ J=a^{-p}\tau.
\]

We first show that images of orthogonal finite-trace projections are
orthogonal and have the same operator norm. The next five lemmas then
construct $u$ and compatible maps $J_e$ on the corners $e\M e$, and prove
that $u$ is unitary. We subsequently combine the maps $J_e$ into a map
on $\mathcal F$, extend it to $\M$, and verify the trace identity above.

\begin{lemma}\anchor{lem:zero-products}
Let $e,f$ be nonzero finite-trace projections with $ef=0$. Then we have $T(e)^*T(f)=T(e)T(f)^*=0$, and $\|T(e)\|_\infty=\|T(f)\|_\infty=\|T(e+f)\|_\infty$.
\end{lemma}

\begin{proof}
Let $e,f\in\mathcal P(\M)\setminus\{0\}$ satisfy $ef=0$ and
$\tau(e),\tau(f)<\infty$. For $\theta\in\mathbb R$, put
\[
A=T(e),\qquad B=T(f),\qquad
C_\theta=A+e^{i\theta}B=T(e+e^{i\theta}f).
\]
Since $ef=fe=0$, we have
\[
|e+e^{i\theta}f|^2
=(e+e^{-i\theta}f)(e+e^{i\theta}f)=e+f,
\qquad |e+e^{i\theta}f|=e+f.
\]
Thus $\frac{e}{\tau(e)^{1/p}}$, $\frac{f}{\tau(f)^{1/p}}$, $\frac{e+e^{i\theta}f}{\tau(e+f)^{1/p}}$ in $\mathscr V$.
By Lemma~\anchorref{lem:keep-form} and linearity of $T$, we have $\frac{A}{\tau(e)^{1/p}}$, $\frac{B}{\tau(f)^{1/p}}$, $\frac{C_\theta}{\tau(e+f)^{1/p}}$ in $\mathscr V$.
Hence $A,B,C_\theta\in\mathcal F\setminus\{0\}$ and $|A|=\|A\|_\infty r(A)$, $|B|=\|B\|_\infty r(B)$, $|C_\theta|=\|C_\theta\|_\infty r(C_\theta)$.
For a finite-trace projection $e$, we have $Q(e)=\tau(e)^{1/p}$.
Since $T$ preserves $Q$, we obtain
\begin{align*}
Q(A)^p&=Q(T(e))^p=Q(e)^p =\tau(e),\\
Q(B)^p&=Q(T(f))^p=Q(f)^p=\tau(f).
\end{align*}
Moreover, $|e+e^{i\theta}f|=e+f$. Thus, by
\eqanchorref{eq:star-values} and \eqanchorref{eq:define-Q},
\begin{align*}
Q(C_\theta)^p
 &=Q(T(e+e^{i\theta}f))^p=Q(e+e^{i\theta}f)^p\\
 &=Q(|e+e^{i\theta}f|)^p=Q(e+f)^p
  =Q(A+B)^p.
\end{align*}
By Lemma~\anchorref{lem:constant-norm}, we have
\[
\|C_\theta\|_\infty=\|A+B\|_\infty\qquad(\theta\in\mathbb R).
\]
Since $|A|=\|A\|_\infty r(A)$ and
$Q(A)^p=\|A\|_\infty^p\tau(r(A))=\tau(e)$, and likewise for
$B$ and $C_\theta$, we obtain, for $j=2,4$,
\begin{equation}\eqanchor{eq:power-values}
\begin{aligned}
\|A\|_j^j
 &=\tau(|A|^j)=\|A\|_\infty^j\tau(r(A))
 =Q(A)^p\|A\|_\infty^{j-p}=\tau(e)\|A\|_\infty^{j-p},\\
\|B\|_j^j
 &=\tau(f)\|B\|_\infty^{j-p},\\
\|C_\theta\|_j^j
 &=\tau(e+f)\|A+B\|_\infty^{j-p}.
\end{aligned}
\end{equation}
Since
\begin{equation}\eqanchor{eq:expand-product}
C_\theta^*C_\theta
=A^*A+B^*B+e^{i\theta}A^*B+e^{-i\theta}B^*A,
\end{equation}
squaring \eqanchorref{eq:expand-product} gives
\begin{equation}\eqanchor{eq:square-product}
\begin{aligned}
(C_\theta^*C_\theta)^2
 &=(A^*A+B^*B)^2+A^*BB^*A+B^*AA^*B\\
 &\quad+e^{i\theta}
 \bigl((A^*A+B^*B)A^*B+A^*B(A^*A+B^*B)\bigr)\\
 &\quad+e^{-i\theta}
 \bigl((A^*A+B^*B)B^*A+B^*A(A^*A+B^*B)\bigr)\\
 &\quad+e^{2i\theta}(A^*B)^2+e^{-2i\theta}(B^*A)^2.
\end{aligned}
\end{equation}
Taking traces in \eqanchorref{eq:expand-product} and integrating gives
\begin{equation}\eqanchor{eq:average-square}
\begin{aligned}
\frac1{2\pi}\int_0^{2\pi}\|C_\theta\|_2^2\,d\theta=\tau(A^*A+B^*B)=\|A\|_2^2+\|B\|_2^2.
\end{aligned}
\end{equation}
On the other hand, taking traces in \eqanchorref{eq:square-product}, we obtain
\begin{equation}\eqanchor{eq:average-fourth}
\begin{aligned}
\frac1{2\pi}\int_0^{2\pi}\|C_\theta\|_4^4\,d\theta&=\tau\bigl((A^*A+B^*B)^2\bigr)
 +\tau(A^*BB^*A)+\tau(B^*AA^*B)\\
 &=\|A\|_4^4+\|B\|_4^4+2\|AB^*\|_2^2+2\|A^*B\|_2^2.
\end{aligned}
\end{equation}

Substituting \eqanchorref{eq:power-values} with $j=2$ into
\eqanchorref{eq:average-square} gives
\begin{equation}\eqanchor{eq:second-powers}
\tau(e+f)\|A+B\|_\infty^{2-p}
=\tau(e)\|A\|_\infty^{2-p}+\tau(f)\|B\|_\infty^{2-p},
\end{equation}
with $j=4$ into \eqanchorref{eq:average-fourth} we have
\begin{equation}\eqanchor{eq:fourth-powers}
\begin{aligned}
\tau(e+f)\|A+B\|_\infty^{4-p}
 &=\tau(e)\|A\|_\infty^{4-p}+\tau(f)\|B\|_\infty^{4-p}+2\|A^*B\|_2^2+2\|AB^*\|_2^2.
\end{aligned}
\end{equation}
Since $(4-p)/(2-p)>1$, the function $t\mapsto t^{(4-p)/(2-p)}$
is strictly convex, and from \eqanchorref{eq:second-powers} we have
\begin{equation}\eqanchor{eq:compare-powers}
\begin{aligned}
&\|A+B\|_\infty^{4-p}
 =\bigl(\|A+B\|_\infty^{2-p}\bigr)^{(4-p)/(2-p)}\\
&\quad=\left(
 \frac{\tau(e)}{\tau(e+f)}\|A\|_\infty^{2-p}
 +\frac{\tau(f)}{\tau(e+f)}\|B\|_\infty^{2-p}
 \right)^{(4-p)/(2-p)}\\
&\quad\le\frac{\tau(e)}{\tau(e+f)}
       \bigl(\|A\|_\infty^{2-p}\bigr)^{(4-p)/(2-p)}
 +\frac{\tau(f)}{\tau(e+f)}
       \bigl(\|B\|_\infty^{2-p}\bigr)^{(4-p)/(2-p)}\\
&\quad=\frac{\tau(e)\|A\|_\infty^{4-p}+\tau(f)\|B\|_\infty^{4-p}}
 {\tau(e+f)}.
\end{aligned}
\end{equation}
Combining this with \eqanchorref{eq:fourth-powers} yields $\|A^*B\|_2^2=\|AB^*\|_2^2=0$. By faithfulness of $\tau$, we obtain
\begin{equation}\eqanchor{eq:zero-products}
A^*B=AB^*=0
\end{equation}
Substituting into
\eqanchorref{eq:fourth-powers} gives
\[
\|A+B\|_\infty^{4-p}
=\frac{\tau(e)}{\tau(e+f)}\|A\|_\infty^{4-p}
 +\frac{\tau(f)}{\tau(e+f)}\|B\|_\infty^{4-p}.
\]
Thus the first and last expressions in \eqanchorref{eq:compare-powers}
are equal, so
\begin{equation}\eqanchor{eq:equality-case}
\begin{aligned}
&\left(
\frac{\tau(e)}{\tau(e+f)}\|A\|_\infty^{2-p}
+\frac{\tau(f)}{\tau(e+f)}\|B\|_\infty^{2-p}
\right)^{\frac{4-p}{2-p}}
=\frac{\tau(e)\|A\|_\infty^{4-p}+\tau(f)\|B\|_\infty^{4-p}}{\tau(e+f)}
\end{aligned}
\end{equation}
If $\|A\|_\infty^{2-p}\ne\|B\|_\infty^{2-p}$, strict convexity gives
\begin{align*}
&\left(
\frac{\tau(e)}{\tau(e+f)}\|A\|_\infty^{2-p}
+\frac{\tau(f)}{\tau(e+f)}\|B\|_\infty^{2-p}
\right)^{\frac{4-p}{2-p}}
<\frac{\tau(e)\|A\|_\infty^{4-p}+\tau(f)\|B\|_\infty^{4-p}}{\tau(e+f)},
\end{align*}
this leads a contradiction. Hence
\[
\|A\|_\infty^{2-p}=\|B\|_\infty^{2-p}.
\]
Substituting this equality into \eqanchorref{eq:second-powers} gives
\begin{equation}\eqanchor{eq:equal-norms}
\|A\|_\infty=\|B\|_\infty=\|A+B\|_\infty.
\end{equation}
\end{proof}

The preceding lemma applies to $e$ and $f-e$ when $e<f$. This gives
the comparison between nested projections used in the next lemma.

\begin{lemma}\anchor{lem:norm-and-support}
For any nonzero finite-trace projections $e,f$, we have
\[
\|T(e)\|_\infty=\|T(f)\|_\infty.
\]
If $0<e\le f$ and $\tau(f)<\infty$, then $r(T(f))=r(T(e))+r(T(f-e))$ and $\ell(T(f))=\ell(T(e))+\ell(T(f-e))$. In particular, $r(T(e))\le r(T(f))$ and $\ell(T(e))\le\ell(T(f))$.
\end{lemma}

\begin{proof}
If $0<e<f$ and $\tau(f)<\infty$, applying \eqanchorref{eq:equal-norms} to $e$ and $f-e$ gives
\begin{equation}\eqanchor{eq:nested-norms}
\|T(e)\|_\infty
=\|T(f-e)\|_\infty
=\|T(e)+T(f-e)\|_\infty
=\|T(f)\|_\infty.
\end{equation}
For arbitrary nonzero finite-trace projections $e,f$,
applying \eqanchorref{eq:nested-norms} to $e\le e\vee f$ and $f\le e\vee f$ gives
$
\|T(e)\|_\infty=\|T(e\vee f)\|_\infty=\|T(f)\|_\infty.
$
For $0<e\le f$ with $\tau(f)<\infty$, Lemma~\anchorref{lem:keep-form} gives
\[
T(e)^*T(e)=\|T(e)\|_\infty^2r(T(e)),\qquad
T(e)T(e)^*=\|T(e)\|_\infty^2\ell(T(e)).
\]
By linearity and \eqanchorref{eq:zero-products} applied to $e$ and $f-e$, we have
\[
T(f)=T(e)+T(f-e),\qquad
T(e)^*T(f-e)=T(e)T(f-e)^*=0.
\]
For $e<f$, \eqanchorref{eq:nested-norms} gives
$\|T(e)\|_\infty=\|T(f-e)\|_\infty=\|T(f)\|_\infty$. Expanding gives
\begin{align*}
\|T(f)\|_\infty^2r(T(f))
 &=T(f)^*T(f)\\
 &=(T(e)+T(f-e))^*(T(e)+T(f-e))\\
 &=T(e)^*T(e)+T(e)^*T(f-e)\\
 &\quad+T(f-e)^*T(e)+T(f-e)^*T(f-e)\\
 &=T(e)^*T(e)+T(f-e)^*T(f-e)\\
 &=\|T(f)\|_\infty^2
   \bigl(r(T(e))+r(T(f-e))\bigr).
\end{align*}
Similarly,
\begin{align*}
\|T(f)\|_\infty^2\ell(T(f))=\|T(f)\|_\infty^2
   \bigl(\ell(T(e))+\ell(T(f-e))\bigr).
\end{align*}
Since $T(f)\ne0$, we have
\begin{align*}
r(T(f))&=r(T(e))+r(T(f-e)),\\
\ell(T(f))&=\ell(T(e))+\ell(T(f-e)).
\end{align*}
\end{proof}

We now use the increasing supports to show that 
$T(e)/\|T(e)\|_\infty$ converge strongly. Their limit will supply the
operator $u$ in the representation of $T$.

\begin{lemma}\anchor{lem:construct-u}
As $e$ runs over the nonzero finite-trace projections, the strong limits
\[
u=\mathop{\mathrm{s\!-\!lim}}_e\frac{T(e)}{\|T(e)\|_\infty},\qquad
u^*=\mathop{\mathrm{s\!-\!lim}}_e\frac{T(e)^*}{\|T(e)\|_\infty}
\]
exist in $\M$. The operator $u$ is a partial isometry satisfying
\[
u^*u=\sup_e r(T(e)),\qquad uu^*=\sup_e\ell(T(e)),
\]
and for every nonzero finite-trace projection $e$, we obtain
$
ur(T(e))=\frac{T(e)}{\|T(e)\|_\infty}.
$
\end{lemma}

\begin{proof}
Let $0<e<f$ with $\tau(f)<\infty$.
Then $(f-e)/\tau(f-e)^{1/p}\in\mathscr V$.
By Lemma~\anchorref{lem:keep-form} and linearity of $T$,
$
\frac{T(f-e)}{\tau(f-e)^{1/p}}\in\mathscr V.
$
Hence
\begin{align*}
T(f-e)^*T(f-e)
 &=\|T(f-e)\|_\infty^2r(T(f-e)),\\
T(f-e)T(f-e)^*
 &=\|T(f-e)\|_\infty^2\ell(T(f-e)).
\end{align*}
Since $\|T(f-e)\|_\infty=\|T(f)\|_\infty$ from
\eqanchorref{eq:nested-norms}, for every vector $\xi$ we obtain
\begin{equation}\eqanchor{eq:vector-difference}
\begin{aligned}
\left\|\frac{(T(f)-T(e))\xi}{\|T(f)\|_\infty}\right\|^2
 &=\frac{\|T(f-e)\xi\|^2}{\|T(f)\|_\infty^2}\\
 &=\frac{\langle T(f-e)^*T(f-e)\xi,\xi\rangle}{\|T(f)\|_\infty^2}\\
 &=\frac{\|T(f-e)\|_\infty^2}{\|T(f)\|_\infty^2}
   \langle r(T(f-e))\xi,\xi\rangle\\
 &=\langle r(T(f-e))\xi,\xi\rangle\\
 &=\langle(r(T(f))-r(T(e)))\xi,\xi\rangle,\\
\left\|\frac{(T(f)^*-T(e)^*)\xi}{\|T(f)\|_\infty}\right\|^2
 &=\frac{\|T(f-e)^*\xi\|^2}{\|T(f)\|_\infty^2}\\
 &=\frac{\langle T(f-e)T(f-e)^*\xi,\xi\rangle}{\|T(f)\|_\infty^2}\\
 &=\frac{\|T(f-e)\|_\infty^2}{\|T(f)\|_\infty^2}
   \langle\ell(T(f-e))\xi,\xi\rangle\\
 &=\langle\ell(T(f-e))\xi,\xi\rangle\\
 &=\langle(\ell(T(f))-\ell(T(e)))\xi,\xi\rangle.
\end{aligned}
\end{equation}
In the suprema below, $g$ ranges over the nonzero finite-trace projections.
Set
\[
q=\sup_g r(T(g)).
\]
By Lemma~\anchorref{lem:norm-and-support}, both support projection nets
are increasing. The standard fact that an increasing net of projections
converges strongly to its supremum gives
\[
\mathop{\mathrm{s\!-\!lim}}_e r(T(e))=q,\qquad
\mathop{\mathrm{s\!-\!lim}}_e\ell(T(e))=\sup_g\ell(T(g)).
\]
For nonzero finite-trace projections $e,f$, take $h=e\vee f$.
The operator norms of $T(e),T(f),T(h)$ are equal. For $k=e,f$,
\eqanchorref{eq:vector-difference} and $r(T(h))\le q$ give
\begin{align*}
\left\|\left(\frac{T(h)}{\|T(h)\|_\infty}
             -\frac{T(k)}{\|T(k)\|_\infty}\right)\xi\right\|^2
 &=\langle(r(T(h))-r(T(k)))\xi,\xi\rangle\\
 &\le\langle(q-r(T(k)))\xi,\xi\rangle.
\end{align*}
The triangle inequality and the preceding estimates now give
\begin{align*}
\left\|\left(\frac{T(e)}{\|T(e)\|_\infty}
             -\frac{T(f)}{\|T(f)\|_\infty}\right)\xi\right\|
 &\le\left\|\left(\frac{T(e)}{\|T(e)\|_\infty}
                  -\frac{T(h)}{\|T(h)\|_\infty}\right)\xi\right\|\\
 &\quad+\left\|\left(\frac{T(h)}{\|T(h)\|_\infty}
                      -\frac{T(f)}{\|T(f)\|_\infty}\right)\xi\right\|\\
 &\le\langle(q-r(T(e)))\xi,\xi\rangle^{1/2}\\
 &\quad+\langle(q-r(T(f)))\xi,\xi\rangle^{1/2}
 \longrightarrow0.
\end{align*}
Here both terms tend to zero because $r(T(e))$ and $r(T(f))$
converge strongly to $q$ as $e,f$ increase.
The same argument with the left supports applies to the adjoints. we have $u^*=\mathop{\mathrm{s\!-\!lim}}_e\frac{T(e)^*}{\|T(e)\|_\infty}$, and 
$u=\mathop{\mathrm{s\!-\!lim}}_e\frac{T(e)}{\|T(e)\|_\infty}$,
where $e$ ranges over the nonzero finite-trace projections.
Products of uniformly bounded strongly convergent nets converge
strongly, so
\[
u^*u=\sup_e r(T(e)),\qquad uu^*=\sup_e\ell(T(e)).
\]
Since $u^*u$ is a projection, $u$ is a partial isometry.
Fix $e$ and let $f\ge e$. Since $e$ and $f-e$ are orthogonal,
$T(f-e)T(e)^*=0$, we obtain
\[
T(f-e)r(T(e))
 =\frac{T(f-e)T(e)^*T(e)}{\|T(e)\|_\infty^2}=0,
\]
and
\begin{align*}
T(f)r(T(e))
 &=\bigl(T(e)+T(f-e)\bigr)r(T(e))
 =T(e).
\end{align*}
The norms $\|T(f)\|_\infty$ and $\|T(e)\|_\infty$ are equal, so we have
\[
\frac{T(f)}{\|T(f)\|_\infty}r(T(e))
 =\frac{T(e)}{\|T(f)\|_\infty}
 =\frac{T(e)}{\|T(e)\|_\infty}.
\]
Taking strong limits over $f\ge e$, using strong continuity of right
multiplication by the fixed projection $r(T(e))$, gives
$
u\,r(T(e))=\frac{T(e)}{\|T(e)\|_\infty}.
$
\end{proof}

The next lemma first constructs $J_e$ on the self-adjoint elements of
$e\M e$ by approximation with finite spectral sums. It then extends
$J_e$ to all of $e\M e$ and checks that choosing a larger finite-trace
corner does not change its value.

\begin{lemma}\anchor{lem:local-map}
For each nonzero finite-trace projection $e$, the formula
\[
J_e(x):=\frac{T(e)^*T(x)}{\|T(e)\|_\infty^2}
\qquad(x\in e\M e)
\]
defines a map $J_e:e\M e\to r(T(e))\M r(T(e))$.
This map is a unital Jordan $*$-monomorphism and satisfies $T(x)=T(e)J_e(x)$ and $\|J_e(b)\|_\infty=\|b\|_\infty$ for $x\in e\M e$ and $ b=b^*\in e\M e$, respectively.
These maps agree on nested corners: if $0<e\le f$ and $\tau(f)<\infty$,
then $J_f(x)=J_e(x)$ for every $x\in e\M e$.
\end{lemma}

\begin{proof}
\textit{Step 1. Self-adjoint elements of $e\M e$.}\\
Fix $e$ and let $0<r\le e$.
The orthogonality of $r$ and $e-r$, from lemmas~\anchorref{lem:zero-products}
and~\anchorref{lem:norm-and-support}, gives $T(e-r)^*T(r)=0$.
So we have
\begin{align*}
\frac{T(e)^*T(r)}{\|T(e)\|_\infty^2}
 &=\frac{\bigl(T(r)+T(e-r)\bigr)^*T(r)}{\|T(e)\|_\infty^2}
 =\frac{T(r)^*T(r)+T(e-r)^*T(r)}{\|T(e)\|_\infty^2}\\
 &=\frac{T(r)^*T(r)}{\|T(e)\|_\infty^2}
 =\frac{\|T(r)\|_\infty^2}{\|T(e)\|_\infty^2}\,r(T(r))
 =r(T(r)).
\end{align*}
The same identity holds for $r=0$.
For a real spectral sum $b=\sum_j\beta_j r_j$ in $e\M e$ with orthogonal nonzero $r_j$, we have
\begin{align*}
T(e)r(T(r_j))
 &=\bigl(T(r_j)+T(e-r_j)\bigr)r(T(r_j))\\
 &=T(r_j)r(T(r_j))+T(e-r_j)r(T(r_j))\\
 &=T(r_j),
\end{align*}
and
\begin{align*}
T(b)
 &=T\left(\sum_j\beta_j r_j\right)
 =\sum_j\beta_j T(r_j)=\sum_j\beta_j T(e)r(T(r_j))
 =T(e)\sum_j\beta_j r(T(r_j)).
\end{align*}

Let $b=b^*\in e\M e$, and choose finite spectral sums $b_n$ of $b$
with $\|b_n-b\|_\infty\to0$. Set
$J_e(b_n)=T(e)^*T(b_n)/\|T(e)\|_\infty^2$. For $n,m$, we write
\[
b_n-b_m=\sum_j\delta_jp_j,\qquad
0\ne p_j\le e,\qquad p_jp_k=0\quad(j\ne k).
\]
Then
\begin{align*}
\|J_e(b_n)-J_e(b_m)\|_\infty
 &=\biggl\|\frac{T(e)^*T(b_n-b_m)}{\|T(e)\|_\infty^2}\biggr\|_\infty\\
 &=\biggl\|\sum\nolimits_j\delta_j\frac{T(e)^*T(p_j)}{\|T(e)\|^2_{\infty}}\biggr\|_\infty\\
 &=\biggl\|\sum\nolimits_j\delta_jr(T(p_j))\biggr\|_\infty\\
 &=\max_j|\delta_j|=\|b_n-b_m\|_\infty\\
 &\rightarrow0.
\end{align*}
Define $J_e(b)$ as its operator-norm limit in $r(T(e))\M r(T(e))$.
Then $T(b_n)=T(e)J_e(b_n)\to T(e)J_e(b)$ in operator norm and also for $Q$, since $Q(y)\le\tau(r(T(e)))^{1/p}\|y\|_\infty$ for $r(y)\le r(T(e))$ and
$Q(T(b_n)-T(b))=Q(b_n-b)\le\tau(e)^{1/p}\|b_n-b\|_\infty\to0$.
Hence
\[
T(b)=T(e)J_e(b),\qquad J_e(b)=\frac{T(e)^*T(b)}{\|T(e)\|_\infty^2}.
\]
The norm equality for finite spectral sums also passes to the limit:
\[
\|J_e(b)\|_\infty
=\lim_n\|J_e(b_n)\|_\infty
=\lim_n\|b_n\|_\infty
=\|b\|_\infty.
\]
We first check that $J_e$ preserves squares of real spectral sums.
For $b=\sum_j\beta_jr_j$, the projections $r_j$ are mutually orthogonal. Thus
\[
b^2=\sum_j\beta_j^2r_j,\qquad
r(T(r_i))r(T(r_j))=0\quad(i\ne j).
\]
Substituting the formula for $b^2$ into the definition of $J_e$ and using
the linearity of $T$ gives the following equalities:
\begin{align*}
J_e(b^2)
 &=\frac{T(e)^*T(b^2)}{\|T(e)\|_\infty^2}
 =\frac{T(e)^*T\bigl(\sum\nolimits_j\beta_j^2r_j\bigr)}{\|T(e)\|_\infty^2}\displaybreak[1]\\
 &=\sum_j\beta_j^2\frac{T(e)^*T(r_j)}{\|T(e)\|_\infty^2}
 =\sum_j\beta_j^2r(T(r_j))\\
 &=\biggl(\sum\nolimits_j\beta_jr(T(r_j))\biggr)^2
 =J_e(b)^2.
\end{align*}
We next pass to any $b=b^*\in e\M e$. Take finite spectral sums $b_n$ of $b$
with $\|b_n-b\|_\infty\to0$. Their squares also converge to $b^2$, since
\[
\|b_n^2-b^2\|_\infty
\le(\|b_n\|_\infty+\|b\|_\infty)\|b_n-b\|_\infty\longrightarrow0.
\]
Each $J_e(b_n)$ is self-adjoint because it is a real linear combination
of projections. Taking limits in operator norm therefore gives
\begin{align*}
J_e(b^2)
 &=\lim_n J_e(b_n^2)
 =\lim_n J_e(b_n)^2
 =J_e(b)^2,\\
J_e(b)^*&=\lim_n J_e(b_n)^*=\lim_n J_e(b_n)=J_e(b).
\end{align*}
\smallskip\noindent\textit{Step 2. All elements of $e\M e$.}\\
Every $x\in e\M e$ has a unique decomposition $x=b+ic$, where
$b=(x+x^*)/2$ and $c=(x-x^*)/(2i)$ are self-adjoint. Define
$J_e(x)=J_e(b)+iJ_e(c)$. The formula from Step~1 and the linearity of $T$ give
\[                                
J_e(x)=\frac{T(e)^*T(x)}{\|T(e)\|_\infty^2},\qquad
T(x)=T(e)J_e(x)\qquad(x\in e\M e).
\]
This formula makes $J_e$ complex-linear. The standard polarization
identities yield a Jordan $*$-map from the
square and self-adjointness identities above; see
\cite[Section~2.2]{FGHS}. In particular,
\[
J_e(x^2)=J_e(x)^2,\qquad J_e(x^*)=J_e(x)^*\qquad(x\in e\M e).
\]
The units of $e\M e$ and $r(T(e))\M r(T(e))$ are $e$ and $r(T(e))$,
respectively, and
\[
J_e(e)=\frac{T(e)^*T(e)}{\|T(e)\|_\infty^2}=r(T(e)),
\]
so $J_e$ preserves the unit. Also, from the identity $T(x)=T(e)J_e(x)$
and the injectivity of $T$ we hold $J_e$ is injective.
Hence $J_e:e\M e\to r(T(e))\M r(T(e))$ is a unital Jordan $*$-monomorphism.
Finally, we check that enlarging the finite-trace corner does not change
the value of the map. For nonzero finite-trace projections $e\le f$,
the identities proved above give
\[
\|T(f)\|_\infty=\|T(e)\|_\infty,\qquad
\frac{T(f)^*T(e)}{\|T(e)\|_\infty^2}=r(T(e)).
\]
For $x\in e\M e$, substitute $T(x)=T(e)J_e(x)$ into the formula for $J_f(x)$:
\begin{align*}
J_f(x)
 &=\frac{T(f)^*T(x)}{\|T(f)\|_\infty^2}\\
 &=\frac{T(f)^*T(e)J_e(x)}{\|T(e)\|_\infty^2}\\
 &=r(T(e))J_e(x)=J_e(x).
\end{align*}
The last equality uses $J_e(x)\in r(T(e))\M r(T(e))$.
Thus $J_f$ agrees with $J_e$ on $e\M e$.
\end{proof}

It remains to show that the partial isometry from
Lemma~\anchorref{lem:construct-u} is unitary. The formulas for
$J_e$ first show that $T$ preserves $\mathcal F$; surjectivity then
forces both support projections of $u$ to be $1$.

\begin{lemma}\anchor{lem:u-is-unitary}
The map $T$ satisfies $T(\mathcal F)=\mathcal F$, and the partial isometry
$u$ of Lemma~\anchorref{lem:construct-u} is unitary.
\end{lemma}

\begin{proof}
Take $x\in\mathcal F$ and a finite-trace projection $e$ such that $x=exe$.
By Lemma~\anchorref{lem:local-map},
\[
T(x)=T(e)J_e(x),\qquad J_e(x)\in r(T(e))\M r(T(e)).
\]
Thus $T(x)$ is bounded, and
\[
r(T(x))\le r(T(e))\le u^*u,\qquad
\ell(T(x))\le\ell(T(e))\le uu^*.
\]
Since $\tau(r(T(e)))<\infty$, we have $T(x)\in\mathcal F$.
Applying Lemmas~\anchorref{lem:zero-products}--\anchorref{lem:local-map}
to $T^{-1}$ gives
$T^{-1}(\mathcal F)\subset\mathcal F$, whence
$\mathcal F\subset T(\mathcal F)$. Therefore $T(\mathcal F)=\mathcal F$.
Every finite-trace projection $q$ is $T(x)$ for some $x\in\mathcal F$,
so $q=r(q)\le u^*u$ and $q=\ell(q)\le uu^*$.
Taking the supremum over such $q$ gives $u^*u=uu^*=1$ by semifiniteness.
Consequently $u$ is unitary.
\end{proof}

Writing $q_e=r(T(e))$, the extensions constructed below fit into the
following commutative diagram:
\[
\begin{tikzcd}[column sep=1.5em, row sep=2em]
e\M e \arrow[r, hook] \arrow[d, "J_e"']
 & \mathcal F \arrow[r, hook] \arrow[d, "S"']
 & \M \arrow[r, hook] \arrow[d, "J"']
 & S(\M,\tau) \arrow[d, "\widetilde J"']
 & \weak \arrow[l, hook'] \arrow[d, "\widetilde J"] \\
q_e\M q_e \arrow[r, hook]
 & \mathcal F \arrow[r, hook]
 & \M \arrow[r, hook]
 & S(\M,\tau)
 & \weak \arrow[l, hook']
\end{tikzcd}
\]
The last two columns describe the measure-continuous extension of $J$
to $S(\M,\tau)$ and its restriction to $\weak$. We use this intermediate
space because $\M\not\subseteq\weak$ when $\tau(1)=\infty$.

We now extend the maps $J_e$ from the finite-trace corners to all of
$\M$ in two steps: combine them into a map $S$ on $\mathcal F$, then
extend $S$ to a normal Jordan $*$-automorphism $J$ of $\M$ and determine
its trace scaling.

\smallskip\noindent\textit{Step 1. Extend to $\mathcal F$.}
Let $a=\|T(e)\|_\infty>0$, independent of the nonzero finite-trace
projection $e$, and set
\[
S(x):=a^{-1}u^*T(x)\qquad(x\in\mathcal F).
\]
For $x\in e\M e$, Lemmas~\anchorref{lem:construct-u}
and~\anchorref{lem:local-map} give
\[
S(x)=a^{-1}u^*T(e)J_e(x)=r(T(e))J_e(x)=J_e(x).
\]
Thus $S$ combines the maps $J_e$ and preserves squares and adjoints.
It is bijective by Lemma~\anchorref{lem:u-is-unitary}. For
$q_e=S(e)=r(T(e))$, the Jordan identity gives
\begin{equation}\eqanchor{eq:compress-S}
S(exe)=q_eS(x)q_e\qquad(x\in\mathcal F).
\end{equation}

For $y\in q_e\M q_e$, put $x=S^{-1}(y)$. Then
\[
S(exe)=q_eS(x)q_e=q_eyq_e=y=S(x).
\]
Hence $x=exe$, proving
\begin{equation}\eqanchor{eq:local-onto}
S(e\M e)=q_e\M q_e.
\end{equation}
Since $S$ and $S^{-1}$ preserve projections, the $q_e$ run over all
finite-trace projections. Moreover, $e\le f$ implies $q_e\le q_f$,
so $e\uparrow1$ and $q_e\uparrow1$ strongly. Lemma~\anchorref{lem:local-map}
also gives
\[
\|S(b)\|_\infty=\|b\|_\infty
\qquad(b=b^*\in\mathcal F).
\]

For $b\in\mathcal F_+$, we have $S(b)=S(b^{1/2})^2\ge0$, and the
same argument applies to $S^{-1}$. Thus positivity of $S=a^{-1}u^*T$
and its inverse is deduced without assuming that $T$ is positive.

\smallskip\noindent\textit{Step 2. Extend to $\M$.}
For $x=x^*\in\M$, define $J(x)$ by its compressions:
\begin{equation}\eqanchor{eq:define-J}
q_eJ(x)q_e=S(exe)\qquad(\tau(e)<\infty).
\end{equation}
Indeed, $\|S(fxf)\|_\infty=\|fxf\|_\infty\le\|x\|_\infty$, so
the net has a weak operator convergent subnet. Its limit satisfies
\eqanchorref{eq:define-J}, since
\[
q_eS(fxf)q_e=S(e(fxf)e)=S(exe)\qquad(f\ge e).
\]
The limit is unique because $q_e\uparrow1$ strongly. The defining
equations then show that $J$ is real-linear on $\M_{\rm sa}$,
$J(1)=1$, and $J(b)=S(b)$ for $b=b^*\in\mathcal F$.

By \eqanchorref{eq:local-onto}, the same construction applies to $S^{-1}$
and gives $K:\M_{\rm sa}\to\M_{\rm sa}$ with
\[
eK(y)e=S^{-1}(q_eyq_e),\qquad K(1)=1.
\]
The compression identities give
\begin{align*}
eK(J(x))e
 &=S^{-1}(q_eJ(x)q_e)=S^{-1}(S(exe))=exe,\\
q_eJ(K(y))q_e
 &=S(eK(y)e)=S(S^{-1}(q_eyq_e))=q_eyq_e.
\end{align*}
Thus $KJ=JK=\mathrm{id}$.

For $x\ge0$, the inequalities $q_eJ(x)q_e=S(exe)\ge0$ imply $J(x)\ge0$;
likewise $K$ is positive. Extend both maps complex-linearly to $\M$.
They are inverse unital positive maps, so Kadison's theorem
\cite{Kadison} makes $J$ a Jordan $*$-automorphism. It extends $S$, hence
\begin{equation}\eqanchor{eq:J-extends-S}
J(x)=S(x),\qquad T(x)=auJ(x)\qquad(x\in\mathcal F).
\end{equation}

It is automatically normal \cite[Section~2.2]{FGHS}.

Finally, for a nonzero finite-trace projection $e$, the isometry of $T$ gives
\[
a\,\tau(J(e))^{1/p}
=Q(auJ(e))=Q(T(e))=Q(e)=\tau(e)^{1/p}.
\]
Thus $\tau(J(e))=a^{-p}\tau(e)$. For any projection $q$, normality gives
\begin{align*}
\tau(J(q))
 =\sup_{f\le q,\,\tau(f)<\infty}\tau(J(f))
 =a^{-p}\sup_{f\le q,\,\tau(f)<\infty}\tau(f)
 =a^{-p}\tau(q).
\end{align*}
Approximation by nonnegative finite spectral sums now yields
\[
\tau(J(x))=a^{-p}\tau(x)\qquad(x\in\M_+),
\]
with infinite values allowed. We have therefore obtained a normal Jordan
$*$-automorphism $J$ with this trace scaling and $T(x)=auJ(x)$ on
$\mathcal F$. It remains to prove
\[
T(x)=au\widetilde J(x)\qquad(x\in\weak).
\]
We first establish the following upper estimate for $Q(b+z_n)$, and use
it in Lemma~\anchorref{lem:equal-distances} to prove that if $Q(x+b)=Q(y+b)$ for every $b\in\mathcal F$, we have $x=y$. 
We then apply this implication to the isometry
$W=(au\widetilde J)^{-1}T$. Since $W$ fixes every element of
$\mathcal F$, it will follow that $W(x)=x$ for all $x\in\weak$,
which proves the required formula.

\begin{lemma}\anchor{lem:limit-bound}
Let $b\in\mathcal F$ and let $(z_n)\subset\weak$ satisfy
$\sup_nQ(z_n)<\infty$ and $z_n\to0$ in measure. Then
\begin{equation}\eqanchor{eq:limit-bound}
\limsup_n Q(b+z_n)\le\max\{Q(b),\limsup_n Q(z_n)\}.
\end{equation}
\end{lemma}

\begin{proof}
If $b=0$, the inequality is immediate. Suppose $b\ne0$ and choose
\[
0<\delta<\epsilon<R,\qquad R>\tau(r(b)).
\]
We estimate $s^{1/p}\mu_s(b+z_n)$ on three intervals.

For $0<s<\epsilon$, boundedness of $b$ gives
\begin{align*}
s^{1/p}\mu_s(b+z_n)
 &\le s^{1/p}\mu_s(z_n)+s^{1/p}\|b\|_\infty\\
 &\le Q(z_n)+\epsilon^{1/p}\|b\|_\infty.
\end{align*}

For $\epsilon\le s\le R$, apply \eqanchorref{eq:sum} we have
\begin{align*}
s^{1/p}\mu_s(b+z_n)
 &\le s^{1/p}\mu_{s-\delta}(b)+s^{1/p}\mu_\delta(z_n)\\
 &\le\left(\frac{s}{s-\delta}\right)^{1/p}Q(b)
       +s^{1/p}\mu_\delta(z_n)\\
 &\le\left(\frac{\epsilon}{\epsilon-\delta}\right)^{1/p}Q(b)
       +R^{1/p}\mu_\delta(z_n).
\end{align*}
For fixed $\delta>0$, convergence in measure gives
$\mu_\delta(z_n)\to0$.

For $s>R$, the finite trace of the support of $b$ gives
\[
\mu_{\tau(r(b))}(b)
\le\|b(1-r(b))\|_\infty=0.
\]
Using \eqanchorref{eq:sum} once more, we obtain
\begin{align*}
\mu_s(b+z_n)
 &\le\mu_{\tau(r(b))}(b)+\mu_{s-\tau(r(b))}(z_n)
 =\mu_{s-\tau(r(b))}(z_n),\\
s^{1/p}\mu_s(b+z_n)
 &\le\left(\frac{s}{s-\tau(r(b))}\right)^{1/p}Q(z_n)
 \le\left(\frac{R}{R-\tau(r(b))}\right)^{1/p}Q(z_n).
\end{align*}

Each bound holds uniformly in $s$ on its interval. Taking the supremum
over $s>0$, followed by the upper limit in $n$, therefore gives
\[
\limsup_n Q(b+z_n)\le
\max\left\{
\begin{aligned}
&\limsup_n Q(z_n)+\epsilon^{1/p}\|b\|_\infty,\\
&\left(\frac{\epsilon}{\epsilon-\delta}\right)^{1/p}Q(b),\\
&\left(\frac{R}{R-\tau(r(b))}\right)^{1/p}\limsup_n Q(z_n)
\end{aligned}
\right\}.
\]
First let $\delta\downarrow0$ with $\epsilon,R$ fixed, then let
$\epsilon\downarrow0$ and $R\uparrow\infty$. This proves
\eqanchorref{eq:limit-bound}.
\end{proof}

\begin{lemma}\anchor{lem:equal-distances}
If $x,y\in \weak$ and $Q(x+b)=Q(y+b)$ for every $b\in\mathcal F$, then $x=y$.
\end{lemma}

\begin{proof}
Suppose $x\ne y$. Choose $|\omega|=1$, $c>0$, and a nonzero finite-trace projection $g\leq1_{[c,\infty)}(\RePart(\omega(y-x)))$. This is possible by semifiniteness, after choosing the phase so that the positive spectral projection is nonzero. Set $h=g\omega(y-x)g-cg$. Since $g^2=g$, we have
\[
gh=hg=h.
\]
Moreover,
\begin{align*}
\RePart h
 &=g\,\RePart(\omega(y-x))\,g-cg\ge cg-cg=0.
\end{align*}
For $N>0$, using $gh=h$ and $h^*g=h^*$ gives
\begin{align*}
\bigl|(N+c)g+h\bigr|^2
 &=\bigl((N+c)g+h^*\bigr)\bigl((N+c)g+h\bigr)\\
 &=(N+c)^2g^2+(N+c)(gh+h^*g)+h^*h\\
 &=(N+c)^2g+(N+c)(h+h^*)+h^*h\\
 &=(N+c)^2g+2(N+c)\RePart h+h^*h\\
 &\ge(N+c)^2g.
\end{align*}
By \eqanchorref{eq:multiply-bound}, $Q(gzg)\le Q(z)$ for $z\in\weak$.
Applying this to $z=Ng+\omega(y-x)$ gives
\begin{align*}
Q(Ng+\omega(y-x))
 &\ge Q\bigl(g(Ng+\omega(y-x))g\bigr)\\
 &=Q\bigl((N+c)g+h\bigr)\\
 &\ge(N+c)\tau(g)^{1/p}.
\end{align*}
Choose $N>0$ with $N\tau(g)^{1/p}>Q(x)$. Set
\[
x_n=x1_{[1/n,n]}(|x|)\in\mathcal F,\qquad z_n=x-x_n.
\]
Then $z_n\to0$ in measure and, by \eqanchorref{eq:multiply-bound}, $Q(z_n)\le Q(x)$.
Now choose the same perturbation for $x$ and $y$:
\[
b_n=\overline\omega Ng-x_n\in\mathcal F,\qquad
\omega(x+b_n)=Ng+\omega z_n.
\]
Since $Q(\omega z_n)=Q(z_n)$, \eqanchorref{eq:limit-bound} yields
\begin{align*}
\limsup_nQ(x+b_n)
 &=\limsup_nQ(Ng+\omega z_n)\\
 &\le\max\{Q(Ng),\limsup_nQ(z_n)\}\\
 &\le\max\{N\tau(g)^{1/p},Q(x)\}
 =N\tau(g)^{1/p}.
\end{align*}
Since $\omega(y+b_n)=Ng+\omega(y-x)+\omega z_n$ converges in measure
to $Ng+\omega(y-x)$, lower semicontinuity of $Q$ \cite[Lemma~3.4]{FK}
and the lower bound above give
\begin{align*}
\liminf_nQ(y+b_n)
 &\ge Q(Ng+\omega(y-x))\\
 &\ge(N+c)\tau(g)^{1/p}>N\tau(g)^{1/p}.
\end{align*}
But $b_n\in\mathcal F$, so the hypothesis gives
$Q(x+b_n)=Q(y+b_n)$ for every $n$. This contradicts the two estimates.
\end{proof}

\begin{proof}[Proof of Theorem~\anchorref{thm:main}]
Fix a nonzero finite-trace projection $e$. By \eqanchorref{eq:scale-Q}, the map
$V=\|T(e)\|_\infty u\widetilde J$ is a surjective isometry.
The same formula proves sufficiency for every triple in the theorem.

The preceding construction gives $T(b)=V(b)$ for $b\in\mathcal F$.
Consequently the surjective linear isometry $W=V^{-1}T$ fixes
$\mathcal F$ pointwise. For $x\in \weak$ and $b\in\mathcal F$,
linearity and isometry give
\[
Q(Wx+b)=Q(Wx+Wb)=Q(W(x+b))=Q(x+b).
\]
Lemma~\anchorref{lem:equal-distances}, applied to $Wx$ and $x$, gives $Wx=x$.
Thus $T=V$ on all of $\weak$; no density of $\mathcal F$ in the
quasi-norm is required.

For uniqueness, suppose a second triple $(a',u',J')$ gives $T$.
For any nonzero finite-trace projection $e$,
$|T(e)|=\|T(e)\|_\infty J(e)=a'J'(e)$. Taking operator norms gives
$a'=\|T(e)\|_\infty$, and hence $J(e)=J'(e)$.
As $e\uparrow1$, normality gives $J(e)\uparrow1$ and $J'(e)\uparrow1$.
The equalities $T(e)/\|T(e)\|_\infty=uJ(e)=u'J'(e)$ then give $u=u'$ by
strong convergence. The two Jordan maps agree on all finite-trace
corners by spectral approximation and linearity. For arbitrary
$x\in\M$, \eqanchorref{eq:compress-J} shows that their
compressions by every $J(e)=J'(e)$ agree. Passing to the strong
limit of these projections proves $J(x)=J'(x)$.
\end{proof}

\section{The finite-trace case}\anchor{sec:finite-trace}

In this section, $0<p<\infty$ and $0<\tau(1)<\infty$.
The von Neumann algebra need not be non-atomic or $\sigma$-finite.
The proof uses the unit $1$, which belongs to $\weak$ in this setting.

\begin{theorem}\anchor{thm:finite}
Let $0<p<\infty$ and let $\M$ be a nonzero von Neumann algebra with a faithful normal trace $\tau$ satisfying $0<\tau(1)<\infty$. A complex-linear map $T\colon\weak\to\weak$ is a surjective isometry for the natural quasi-norm if and only if there are a unitary $u\in\M$ and a trace-preserving normal Jordan $*$-automorphism $J$ of $\M$ such that
\[
T(x)=u\widetilde J(x),\qquad x\in\weak,
\]
where $\widetilde J$ denotes the extension of $J$ to $\tau$-measurable operators. The pair $(u,J)$ is unique, and $u=T(1)$.
\end{theorem}

Normalize $\tau(1)=1$, which multiplies the original quasi-norm by a fixed constant. The supremum over $s$ in \eqanchorref{eq:define-Q} may now be restricted to $0<s<1$, and $Q(1)=1$. The proof recovers the order on $\weak$ from the behavior of $Q$ near this unit.

\begin{lemma}\anchor{lem:real-part}
For $x\in \weak$, the inequality $\RePart x\ge0$ holds if and only if $Q(1+rx)\ge1$ for all sufficiently small $r>0$.
\end{lemma}

\begin{proof}
If $\RePart x\ge0$, then in $S(\M,\tau)$ we have
\[
|1+rx|^2=1+2r\RePart x+r^2x^*x\ge1,\qquad r>0.
\]
It follows by spectral calculus that $|1+rx|\ge1$, so $Q(1+rx)\ge1$.

Suppose $\RePart x\not\ge0$. Choose $c>0$ such that $e=1_{(-\infty,-c]}(\RePart x)\ne0$. Choose $M>0$ such that $q=e\wedge 1_{[0,M]}(|x|)\ne0$. This is possible since the trace of $1_{(M,\infty)}(|x|)$ tends to zero and $\tau(e\wedge f)\ge\tau(e)-\tau(1-f)$ for projections $e,f$. Thus $q\RePart xq\le-cq$ and $\|xq\|_\infty\le M$. Compressing $|1+rx|^2$ gives
\[
\|(1+rx)q\|_\infty^2\le1-2cr+r^2M^2.
\]
For $1-\tau(q)\le s<1$, formula~\eqanchorref{eq:compute-mu} and $s^{1/p}\le1$ give $s^{1/p}\mu_s(1+rx)\le\|(1+rx)q\|_\infty$. For $0<s<1-\tau(q)$, the same formula gives $\mu_s(1+rx)\le1+r\mu_s(x)$, and hence
\[
s^{1/p}\mu_s(1+rx)\le(1-\tau(q))^{1/p}+rQ(x).
\]
Here $M>0$ because $q\RePart xq\le-cq$ and $q\ne0$.
Choose
\[
0<r<\min\left\{\frac{2c}{M^2},\,
\frac{1-(1-\tau(q))^{1/p}}{Q(x)}\right\}.
\]
The denominator $Q(x)$ is positive since $\RePart x$ has a negative
part. The first bound is then at most
$\sqrt{1-2cr+r^2M^2}<1$, and the second is strictly less than one
by the second restriction on $r$. If $\tau(q)=1$, the second range
of scales is empty and can simply be omitted.
Taking the supremum over the two ranges proves $Q(1+rx)<1$.
\end{proof}

\begin{lemma}\anchor{lem:find-unitaries}
An element $y\in \weak$ with $Q(y)=1$ is unitary if and only if for every nonzero $z\in \weak$ and every $\eta>0$ there is a scalar $\lambda$ with $0<|\lambda|<\eta$ such that $Q(y+\lambda z)<1$.
\end{lemma}

\begin{proof}
For $y=1$, choose $\omega\in\{1,-1,i,-i\}$ such that $\RePart(\omega z)\not\ge0$. Such an $\omega$ exists because $z\ne0$. Lemma~\anchorref{lem:real-part} gives the desired perturbation along $\omega z$. Left multiplication by a unitary preserves $Q$, so the property holds for every unitary $y$.

Conversely, write $y=vh$ with $h=|y|$. Since $\tau$ is finite, $v$ extends to a unitary \cite[Lemma~2.5]{dJC}. Thus it suffices to verify the condition for $h$, which also satisfies $Q(h)=1$. Suppose that $h\ne1$.

If $h\not\ge1$, choose $0<b<1$ such that $q=1_{[0,b]}(h)\ne0$. In \eqanchorref{eq:compute-Q} thresholds $a\le b$ contribute at most $b$, and thresholds $a>b$ give the same spectral projections for $h$ and $h(1-q)$. Consequently $Q(h(1-q))=1$. Right multiplication by $1-q$ gives
\[
Q(h+\lambda q)\ge Q(h(1-q))=1,\qquad \lambda\in\mathbb C.
\]
If $h\ge1$, choose $\varepsilon>0$ such that $q=1_{(1+\varepsilon,\infty)}(h)\ne0$. For $|\lambda|<\varepsilon$, spectral calculus gives $|h+\lambda q|\ge1$, so $Q(h+\lambda q)\ge1$. In each case $z=q$ contradicts the perturbation property.
\end{proof}

\begin{proof}[Proof of Theorem~\anchorref{thm:finite}]
Let $T$ be as in Theorem~\anchorref{thm:finite}. By Lemma~\anchorref{lem:find-unitaries},
$u:=T(1)$ is unitary. Define $S(x)=u^*T(x)$; this is a unital
surjective isometry. Lemma~\anchorref{lem:real-part} gives
\[
\RePart S(x)\ge0\quad\Longleftrightarrow\quad\RePart x\ge0.
\]
An element $x$ is self-adjoint exactly when both $\RePart(ix)$ and $\RePart(-ix)$ are nonnegative. Complex linearity and the displayed equivalence therefore show that $S$ preserves self-adjointness and positivity in both directions.
This derives positivity of $S=u^*T$ and its inverse directly from the
quasi-norm, without assuming that $T$ is positive.

Every bounded self-adjoint element lies in some order interval $[-c1,c1]$, and these intervals contain only bounded elements. Since $S$ preserves each such interval in both directions, it restricts to a unital linear order isomorphism on $\M$. By Kadison's theorem \cite{Kadison}, this restriction is a Jordan $*$-automorphism $J$. It is automatically normal \cite[Section~2.2]{FGHS}. For every projection $e\in\M$,
\[
\tau(J(e))^{1/p}=Q(J(e))=Q(e)=\tau(e)^{1/p}.
\]
By normality and spectral approximation, $J$ preserves $\tau$.

Let $\widetilde J$ be the trace-preserving extension of $J$ to
$S(\M,\tau)$ supplied by \cite[Lemma~2.2]{dJC}. Formula
\eqanchorref{eq:scale-Q}, with $\kappa=1$, shows that $\widetilde J$ and
its inverse map $\weak$ into $\weak$. They preserve positivity, and hence
are order isomorphisms on the self-adjoint part of $\weak$.
For $x\in(\weak)_+$ put $x_n=x\wedge n1$. Then $x_n\uparrow x$ in the self-adjoint part of $\weak$.
Since $S$ and $\widetilde J$ are order isomorphisms, they preserve these
suprema. It follows that
\[
S(x)=\sup_n S(x_n)=\sup_n J(x_n)
=\sup_n\widetilde J(x_n)=\widetilde J(x).
\]
For self-adjoint $x\in \weak$, its positive and negative parts belong
to $\weak$, since their singular values are bounded by those of $x$.
Writing $x=x_+-x_-$ extends equality to self-adjoint elements.
For general $x$, \eqanchorref{eq:sum} and \eqanchorref{eq:star-values} put
$\RePart x$ and $\Im x=(x-x^*)/(2i)$ in $\weak$; use
$x=\RePart x+i\Im x$ and complex linearity.
This proves $S=\widetilde J$ on $\weak$.
Finally any representation $T=u\widetilde J$ has $T(1)=u$ because
$J(1)=1$. Thus $u$ is unique, and $J=(u^*T)|_{\M}$ is unique as well.

Conversely, \eqanchorref{eq:scale-Q} with $\kappa=1$ shows that
$x\mapsto u\widetilde J(x)$ is a surjective isometry.
\end{proof}

Thus the finite-trace representation has no scalar factor: if a representation
$T=au\widetilde J$ satisfies $\tau\circ J=a^{-p}\tau$, evaluating at $1$
gives $\tau(1)=a^{-p}\tau(1)$ and hence $a=1$.

\section*{Declaration of AI use}

The organization and ordering of the lemmas, together with the drafting
and development of their proofs, were carried out by GPT-6.1 Sol under
the authors' direction. The tool was also used to refine the exposition,
clarify intermediate calculations, and assist with LaTeX preparation.
The authors take full responsibility for the mathematical correctness
and readability of the lemmas and proofs, the accuracy of the references,
and the final content of the manuscript. All decisions concerning the
statements, arguments, and presentation remain the responsibility of
the authors. The AI tool is not listed as an author.

\end{document}